\documentclass[10pt,a4paper]{article}
\usepackage[top=70pt,bottom=90pt,left=64pt,right=64pt]{geometry}
\usepackage{amsmath} 
\usepackage{amsthm}
\usepackage{amssymb} 
\usepackage{amsfonts}
\usepackage{mathrsfs}
\usepackage{graphicx} 
\usepackage{inputenc}
\usepackage{upref} 
\usepackage{float}
\usepackage[all]{xy}
\usepackage{bm}
\usepackage{calc}
\usepackage{hyperref}
\usepackage{tikz-cd}
\usepackage[all]{xy}
\usepackage{hyperref}

\def\makeautorefname#1#2{\expandafter\def\csname#1autorefname\endcsname{#2}}
\def\equationautorefname~#1\null{(#1)\null}

\makeautorefname{theorem}{Theorem}
\makeautorefname{lemma}{Lemma}
\makeautorefname{proposition}{Proposition}
\makeautorefname{remark}{Remark}
\makeautorefname{definition}{Definition}
\makeautorefname{conjecture}{Conjecture}
\makeautorefname{construction}{Construction}
\makeautorefname{corollary}{Corollary}

\newtheorem{theorem}{Theorem}[section]

\newtheorem{lemma}{Lemma}[section]
\newtheorem{definition}{Definition}[section]

\newtheorem{proposition}{Proposition}[section]
\newtheorem{corollary}{Corollary}[section]
\newtheorem{remark}{Remark}[section]

\makeatletter
\let\c@lemma=\c@theorem
\let\c@proposition=\c@theorem
\let\c@remark=\c@theorem
\let\c@definition=\c@theorem
\let\c@conjecture=\c@theorem
\let\c@construction=\c@theorem
\let\c@corollary=\c@theorem
\numberwithin{equation}{section}
\makeatother

\newcommand{\Z}{\mathbb{Z}}

\title{Universal Spaces and Splittings of Equivariant Spectra}
\author{Yutao Liu \\ University of Chicago
\\ Email: \href{mailto:yutao492@math.uchicago.edu}{yutao492@math.uchicago.edu}
}

\begin{document}
\maketitle

\paragraph{Abstract:} Let $G$ be a finite group. We re-analyze the splitting of rational $G$-spectra, which is discussed in Barne's thesis and traces back to Greenlees and May. We will study how the splitting behaves under a localization which is weaker than rationalization and discuss its application in computing equivariant cohomology of $G$-spectra. In particular, we will explicitly compute $\pi_\bigstar(H\underline{\Z})$ and $\pi_\bigstar(HA_G)$ for the dihedral group $G=D_{2p}$. The additive structures and partial multiplicative structures are already computed by Kriz, Lu, and Zou. Our new method will provide a systematic way to compute them as $RO(D_{2p})$-graded rings.

\tableofcontents

\section{Introduction}

For any finite group $G$, the category of rational $G$-spectra admits a splitting, which is first discussed by Greenlees and May:

\begin{theorem} \label{idemsplit}

\textbf{[GM95]} There is an orthogonal basis $\{e_H:H\subset G\}$ of the rational Burnside ring, which only contains idempotent elements. For any rational $G$-spectrum $X$, define $e_HX$ as
$$e_HX:=colim(X\xrightarrow{e_H}X\xrightarrow{e_H}X\xrightarrow{e_H}...)$$
Then we have
$$X\simeq\bigvee_H e_HX$$
and
$$[X,Y]^G\cong\prod_H[e_HX,e_HY]^G$$
with one $H$ chosen from each conjugacy class of subgroups.
\end{theorem}

This theorem is further studied by Barnes in \textbf{[Bar08]}, which reproves the splitting above by some topological constructions and sets up an algebraic model which completely describes the behavior of rational $G$-spectra.

One important idea in Barnes' proof is to use the universal $G$-space $E\lambda_H$, which is characterized by the fixed point subspaces:
$$(E\lambda_H)^K\simeq
\begin{cases}
S^0,\text{ if } K \text{ is conjugate to } H,\\
*, \text{ otherwise.}
\end{cases}$$

\begin{theorem} \label{computesplit}

\textbf{[Bar08]} For rational $G$-spectra $X,Y$, we have an isomorphism
$$[X,Y]^G\cong\prod_H [E\lambda_H\wedge X,E\lambda_H\wedge Y]^G$$
with one $H$ chosen from each conjugacy class. Moreover, the functor
$$X\mapsto\bigvee_H E\lambda_H\wedge X$$
is symmetric monoidal.

\end{theorem}
\medskip

However, the rational assumption is too strong. In this paper, we will reprove this theorem by further studying the construction and properties of universal spaces. We will show:

\begin{theorem} \label{orderinvert}

The splitting in \autoref{computesplit} happens when we only invert all prime factors of $|G|$, instead of applying a full rationalization.

\end{theorem}
\bigskip

Moreover, if we only invert some, but not all prime factors of $|G|$, we may still get a partial splitting
$$[X,Y]^G\cong[E\mathscr{F}_+\wedge X,E\mathscr{F}_+\wedge Y]^G\oplus [\widetilde{E\mathscr{F}}\wedge X,\widetilde{E\mathscr{F}}\wedge Y]^G$$
where $E\mathscr{F}$ is the universal space corresponding to the family $\mathscr{F}$.

The corresponding functor
$$X\mapsto (E\mathscr{F}_+\wedge X)\vee(\widetilde{E\mathscr{F}}\wedge X)$$
still turns out to be symmetric monoidal.
\medskip

The summands $E\mathscr{F}_+\wedge X$ and $\widetilde{E\mathscr{F}}\wedge X$ are usually more simple than $X$ since the product with either $E\mathscr{F}_+$ or $\widetilde{E\mathscr{F}}$ reduces the types of cells in $X$. When we want to compute the $RO(G)$-graded homotopy $\pi_\bigstar(X)$ of $X$, we can break it into more computable parts:
$$\pi_\bigstar(X)\cong\pi_\bigstar(E\mathscr{F}_+\wedge X)\oplus\pi_\bigstar(\widetilde{E\mathscr{F}}\wedge X)$$
after inverting proper prime factors of $|G|$. Inverting different prime factors in $|G|$ gives us different partial splittings. We can glue these data back to recover the unlocalized $\pi_\bigstar(X)$.
\medskip

One explicit example given in this paper is the computation of $\pi_\bigstar(H\underline{\Z})$ when $G=D_{2p}$ is a dihedral group. The additive structure is already computed in \textbf{[KL20]} and \textbf{[Zou18]}, while the latter also provides some partial multiplicative structure. With our new method, the complete multiplicative structure follows naturally since all functors in the splittings are symmetric monoidal.

The main idea is, we get two different partial splittings when $2$ and $p$ are inverted respectively. Both splittings break $\pi_\bigstar(H\underline{\Z})$ into two more computable parts, whose computations can be decomposed into some $RO(C_2)$ or $RO(C_p)$-graded cohomology theories. We can compute $\pi_\bigstar(H\underline{\Z})[1/2]$ and $\pi_\bigstar(H\underline{\Z})[1/p]$ from these two splittings, and hence recover $\pi_\bigstar(H\underline{\Z})$.

\bigskip

\paragraph{Structure of the paper:}

In section 2, we will introduce the universal space $E\mathscr{F}$. We are especially interested in the equivariant homology of $E\mathscr{F}$ with coefficients in the Burnside ring Mackey functor $A_G$. The most important homological information is given in \autoref{degree0} and \autoref{torsion}.

Section 3 contains a reproof of \autoref{computesplit} and a proof for \autoref{orderinvert}. We will show that the splitting in \autoref{computesplit} can be obtained by repeatedly applying the partial splitting
$$X\simeq (E\mathscr{F}_+\wedge X)\vee(\widetilde{E\mathscr{F}}\wedge X)$$
for all choices of $\mathscr{F}$. Since each partial splitting only requires some prime factors of $|G|$ to be inverted, \autoref{computesplit} works when we invert $|G|$. We will also give a criterion about which prime factors should be inverted in order to make the partial splitting happen. This criterion will be quite useful for our later computations.

In sections 4 and 5, we consider the dihedral group $G=D_{2p}$ and compute the $RO(G)$-graded homotopy of the Eilenberg-Maclane spectra $H\underline{\Z}$. We will discuss how this method can be applied to more general $D_{2p}$-spectra in section 6. In particular, we compute the $RO(G)$-graded homotopy of $HA_G$.

Finally, in section 7 and 8, we will explain how this idea can be applied for more general finite groups. We will focus on the case that $G=G_1\ltimes G_2$, and prove two partial splittings when either $|G_1|$ or $|G_2|$ is inverted.
\bigskip

Here is a list of main results of this paper:

The homotopy of $H\underline{\Z}$ as a $RO(D_{2p})$-graded ring: \autoref{hz}.

The homotopy computation of more general $D_{2p}$-spectra: \autoref{generald2p} and \autoref{hagrog}.

The splitting theorem for more general finite groups: \autoref{bigtwosplit}.
\bigskip

\paragraph{Notations:} In this paper, we use $*$ when the homotopy or homology is graded over $\Z$, and use $\bigstar$ when graded over $RO(G)$. We add an underline to express the Mackey functor valued homotopy or homology: $\underline{H}$, $\underline{\pi}$.

To be more precise, let $X$ be a $G$-space and $M$ be a Mackey functor. Define the Mackey functor valued equivariant homology of $X$ with coefficients in $M$ as
$$\underline{H}_*^G(X;M):=\underline{\pi}_*^G(X\wedge HM).$$
So we have
$$\underline{H}_*^G(X;M)(G/H)=[\Sigma^*G/H_+,HM\wedge X]^G,$$
$$H_*^G(X;M)=\underline{H}_*^G(X;M)(G/G).$$
Here $HM$ is the equivariant Eilenberg-Maclane spectrum corresponding to $M$.
\medskip

We discuss the properties of the universal spaces by computing their $\Z$-graded, Mackey functor valued homology $\underline{H}_*$ in section 2 and 3. But we compute the $RO(G)$-graded, abelian group (to be more precise, graded ring) valued homotopy $\pi_\bigstar$ in the remaining sections. Theoretically, our computational methods also work for $\underline{\pi}_\bigstar$, although the multiplicative structures expressed as box products may become quite complicated.

\section{Universal spaces}

We will introduce the universal space $E\mathscr{F}$ and study its equivariant homology in this section. The most important properties are given in \autoref{homsplit} and \autoref{torsion}.

\begin{definition} \label{clspace}
A \textbf{family} $\mathscr{F}$ is a collection of subgroups of $G$ which is closed under conjugation and taking subgroups. The corresponding \textbf{universal space} $E\mathscr{F}$ is an unbased $G$-space such that $(E\mathscr{F})^H$ is contractible for all $H\in\mathscr{F}$ and empty otherwise.
\end{definition}

\begin{remark}
When $\mathscr{F}$ only contains the trivial subgroup, $E\mathscr{F}$ becomes $EG$, which is a contractible space with free $G$-action.
\end{remark} 
\smallskip

One construction of $E\mathscr{F}$ is given in \textbf{[Dieck72]} as a join of spaces. We will explain this idea later and use it to prove some important properties. An alternative construction is given by \textbf{[Elm83]} as a categorical bar construction.
\smallskip

We can fully characterize $E\mathscr{F}$ by fixed point subspaces:

\begin{lemma} \label{unique}
The universal space $E\mathscr{F}$ is unique up to weak equivalence.
\end{lemma}

\paragraph{Proof:} Assume that $X,Y$ are $G$-spaces with the same conditions on fixed point subspaces as $E\mathscr{F}$.

If we have a $G$-map $X\rightarrow Y$, its restriction on each fixed point subspace is a weak equivalence since $X,Y$ have the same types of fixed point subspaces as either $*$ or $\emptyset$. Thus this map is a $G$-weak equivalence.

If we do not have a map between $X,Y$, consider $X\times Y$, which is another $G$-space with the same conditions on fixed point subspaces. The above argument shows that the projections $X\times Y\rightarrow X$ and $X\times Y\rightarrow Y$ are $G$-weak equivalences. Therefore, $E\mathscr{F}$ is unique up to weak equivalence. $\Box$
\bigskip

We assume $E\mathscr{F}$ to be a $G$-CW complex by applying the $G$-CW approximation.
\begin{lemma} \label{uniqueef}
Any $G$-self map of $E\mathscr{F}$ is $G$-homotopic to the identity.
\end{lemma}

This lemma is implied by the following theorem:

\begin{theorem} \label{elmendorf}
Let $\mathcal{O}_G$ be the category of $G$-orbits. Define an $\mathcal{O}_G$-space to be a contravariant functor from $\mathcal{O}_G$ to topological spaces.

The following pair of functors
$$\Phi: G-spaces\rightleftarrows\mathcal{O}_G-spaces:\Psi$$
defined by $\Phi(X)(G/H):=X^H$ and $\Psi(T):=T(G/\{e\})$ form a Quillen equivalence.
\end{theorem}

We refer to [\textbf{May96}, VI.6] for more details.

\paragraph{Proof of \autoref{uniqueef}:} Construct an $\mathcal{O}_G$-space $T$ by setting 
$$T(G/H)=
\begin{cases}
*, \text{ if }H\in\mathscr{F}\\
\emptyset, \text{ otherwise}
\end{cases}$$

Let $CT$ be a cofibrant approximation of $T$. Then
$$\Psi(CT)^H=\Phi\Psi(CT)(G/H)\cong CT(G/H)\cong
\begin{cases}
*, \text{ if }H\in\mathscr{F}\\
\emptyset, \text{ otherwise}
\end{cases}$$

Recall that
$$\Phi E\mathscr{F}(G/H)=(E\mathscr{F})^H\cong
\begin{cases}
*, \text{ if }H\in\mathscr{F}\\
\emptyset, \text{ otherwise}
\end{cases}$$

The unique map $\Phi E\mathscr{F}\rightarrow T$ induces $\Phi E\mathscr{F}\rightarrow CT$, whose adjuction $E\mathscr{F}\rightarrow\Psi(CT)$ becomes a $G$-equivalence.

So we have
$$[E\mathscr{F},E\mathscr{F}]^G\cong[E\mathscr{F},\Psi(CT)]^G\cong[\Phi E\mathscr{F},CT]^{\mathcal{O}_G}\cong[\Phi E\mathscr{F},T]^{\mathcal{O}_G},$$
which contains a single element.

Thus the identity map is the only self-map of $E\mathscr{F}$ up to $G
$-homotopy. $\Box$

\bigskip

Notice that $E\mathscr{F}\times E\mathscr{F}$ also appears as a universal space for $\mathscr{F}$. Thus the above two lemmas imply:

\begin{lemma} \label{monoid}
The universal space $E\mathscr{F}$ is a $G$-topological semigroup up to homotopy.
\end{lemma}

Usually $E\mathscr{F}$ is not a Hopf space since there is no unit. But with extra conditions, $\Sigma^\infty E\mathscr{F}_+$ will appear as a homotopy ring spectrum. We will explain this idea in section 3.
\medskip

Now we discuss the homological properties for $E\mathscr{F}$. First we recall the definition of the Burnside ring and the corresponding Mackey functor:

\begin{definition} \label{ag}

For any finite group $G$, the collection of isomorphism classes of finite $G$-sets forms a commutative monoid, with addition induced by disjoint unions. The \textbf{Burnside ring} $A(G)$ is defined to be the group completion of this monoid. We can express $A(G)$ as a free $\Z$-module, whose basis corresponds to $G$-orbits. We use $\{G/H\}$ to denote the basis element of the orbit $G/H$.
\medskip

The \textbf{Burnside ring Mackey functor} $A_G$ is defined by $A_G(G/H):=A(H)$. The Mackey functor structure is given by the following maps:
\smallskip

For any $L\subset H$, the transfer map 
$$T_L^H: A(L)\rightarrow A(H)$$
sends $\{L/K\}$ to $\{H/K\}$ for any $K\subset L$. 

The restriction map
$$R_L^H: A(H)\rightarrow A(L)$$
sends $\{H/N\}$ to itself, but viewed as an $L$-space.

When $L=g^{-1}Hg$, we also have an isomorphism
$$C_L^H:A(L)\rightarrow A(H)$$
sending $\{L/K\}$ to $\{H/gKg^{-1}\}$ for any $K\subset L$.
\medskip

The Cartesian product of finite $G$-sets makes $A(G)$ into a ring. This multiplicative structure also makes $A_G$ into a Green functor, which is a monoid under the box product.
\end{definition}

\begin{proposition} \label{burnsidehom}
For a $G$-CW complex $X$, the integer graded homology of $X$ with coefficients in $A_G$ can be computed as
$$H_*^G(X;A_G)\cong\bigoplus_H H_*(X^H/W_GH;\Z)$$
with one $H$ chosen from each conjugacy class of subgroups of $G$. Here $W_GH$ is the Weyl group of $H$ in $G$.

As a Mackey functor,
$$\underline{H}_*^G(X;A_G)(G/L)\cong H_*^L(X;A_L)\cong\bigoplus_K H_*(X^K/W_LK;\Z)$$
with one $K$ chosen from each conjugacy class of subgroups of $L$.
\end{proposition}

\paragraph{Proof:} It suffices to prove the first equation. The second one can be proved by the same argument.
\medskip

When we view $A_G$ as a covariant coefficient system, the restriction maps are removed. For any $G$-map $G/L_1\rightarrow G/L_2$ induced by the multiplication of $g\in G$ (where we require $g^{-1}L_1g\subset L_2$), we have an induced map
$$A(L_1)\rightarrow A(L_2)$$
sending each $\{L_1/H\}$ to $\{L_2/g^{-1}Hg\}$. Therefore, as a coefficient system, $A_G$ can be decomposed as
$$A_G=\bigoplus_{H}A_G^H$$
with one $H\subset G$ chosen from each conjugacy class. Here $A_G^H$ is the sub-coefficient system of $A_G$ such that each $A_G^H(G/L)\subset A(L)$ is generated by all $\{L/g^{-1}Hg\}$ with $g^{-1}Hg\subset L$. 
\smallskip

Now we have
$$H_*^G(X;A_G)\cong\bigoplus_{H}H_*^G(X;A_G^H).$$
It suffices to prove
$$H_*^G(X;A_G^H)\cong H_*(X^H/W_GH;\Z)$$
for any $H\subset G$.
\bigskip

As described in \textbf{[Wil75]}, the equivariant homology can be computed in a cellular way:

Let $\underline{C}_*(X)$ be the chain complex of contravariant coefficient systems such that
$$\underline{C}_*(X)(G/H):=C_*(X^H),$$
where $C_*(X^H)$ is the cellular chain complex of $X^H$ and the orbit maps are sent to conjugacies and inclusions of subcomplexes. For any covariant coefficient system $M$, $H_*^G(X;M)$ is the homology of the chain complex:
$$C_*^G(X;M):=\underline{C}_*(X)\otimes_{\mathcal{O}_G}M,$$
where the tensor product is taken over the category $\mathcal{O}_G$ of $G$-orbits. To be more precise, we have
$$C_*^G(X;M):=\left(\bigoplus_{L\subset G}C_*(X^L)\otimes M(G/L)\right)/\sim.$$

For any $G$-map $f:G/L_1\rightarrow G/L_2$, the equivalence relation identifies $f^*a\otimes b$ and $a\otimes f_*b$ for all $a\in C_*(X^{L_2})$ and $b\in M(G/L_1)$.
\medskip

Now we choose $M=A_G^H$. Since each $M(G/L)$ is freely generated by $\{L/g^{-1}Hg\}$ and the structure maps send these basis elements to each other, $C_*^G(X;M)$ is the free $\Z$-module generated by the equivalence classes of $e\otimes \{L/g^{-1}Hg\}$, for all cells $e$ in $X^L$ and $g^{-1}Hg\subset L$.
\medskip

We can eliminate the equivalence relation by the following four facts:

\textbf{(1)} The free $\Z$-module $A_G^H(G/H)$ is generated by a single element $\{H/H\}$.

\textbf{(2)} For any $e\otimes\{L/g^{-1}Hg\}$, it is identified with $ge\otimes\{H/H\}$ by the equivalence relation.

\textbf{(3)} If there exists another $g^\prime\in G$ such that
$$ge\otimes\{H/H\}\sim e\otimes\{L/g^{-1}Hg\}\sim g^\prime e\otimes\{H/H\},$$
then $\{L/g^{-1}Hg\}=\{L/(g^{\prime})^{-1}Hg^\prime\}$. Thus
$$g^{-1}Hg=l^{-1}((g^\prime)^{-1}Hg^\prime)l$$
for some $l\in L$. So $g^\prime lg^{-1}\in W_GH$. Since $e$ is $L$-fixed, we have $(g^\prime lg^{-1})(ge)=g^\prime e$. Thus $ge$ and $g^\prime e$ are in the same $W_GH$-orbit.

\textbf{(4)} The converse of \textbf{(3)} is also true: $e_1\otimes\{H/H\}$ and $e_2\otimes\{H/H\}$ are identified if $e_1,e_2$ are in the same $W_GH$-orbit.
\medskip

Now we have a 1-1 correspondence between the basis of $C_*^G(X;A_G^H)$ and the cells in $X^H/W_GH$. Therefore, we get
$$C_*^G(X;A_G^H)\cong C_*(X^H/W_GH).$$
Taking the homology on both sides gives us the required equation. $\Box$
\bigskip

\begin{remark}
\autoref{burnsidehom} follows from the fact that the underlying coefficient system of $A_G$ splits into $A_G^H$. However, such splitting cannot be lifted to the Mackey functor level. Thus the decomposition does not work for the cohomology with coefficients in $A_G$.

In fact, since we cannot define orbit spectra on the complete $G$-universe,
$$X\mapsto H_*(X^H/W_GH)$$
is a homology theory only for $G$-spaces, and hence cannot be represented by any $G$-spectrum.
\end{remark}

Since $E\mathscr{F}$ only has empty and contractible fixed point subspaces, \autoref{burnsidehom} helps us to compute the 0th degree equivariant homology of $E\mathscr{F}$ explicitly:

\begin{proposition} \label{degree0}
The trivial map $E\mathscr{F}\rightarrow\{*\}$ implies an inclusion
$$\underline{H}_0^G(E\mathscr{F},A_G)\hookrightarrow\underline{H}_0^G(*,A_G)\cong A_G.$$

The image of $\underline{H}_0^G(E\mathscr{F},A_G)(G/L)$ in $A(L)$ is generated by all $\{L/K\}\in A(L)$ with $K\in\mathscr{F}$. 
\end{proposition}

We already have a small splitting here:

\begin{lemma} \label{homsplit}
When $|G|$ is inverted, $M_\mathscr{F}:=\underline{H}_0^G(E\mathscr{F};A_G)$ is a direct summand of $A_G$ as Mackey functors.
\end{lemma}

\paragraph{Proof:} Assume that $|G|$ is inverted everywhere. For any $H\subset G$, define a linear map
$$\varphi_H:A(G)\rightarrow\Z$$
which sends each $G$-set $S$ to $|S^H|$. 

Let $s_{(K,H)}$ denote $\varphi_K(\{G/H\})=|G/H|^K$, which can be computed as the product between $|W_GH|$ and the number of subgroups of $G$ containing $K$ and in the conjugacy class of $H$. Thus $s_{(K,H)}\neq 0$ if and only if $K$ is sub-conjugate to $H$.

When $|G|$ is inverted, all $|W_GH|$ and non-zero $s_{(K,H)}$ become invertible. We choose elements $e_H\in A(G)$ inductively by defining
$$e_H=|W_GH|^{-1}\left(\{G/H\}-\sum_K s_{(K,H)}^{-1}e_K\right)$$
with one $K$ chosen from each conjugacy class that contains a proper subgroup of $H$.

By induction, we have
$$\varphi_K(e_H)=
\begin{cases}
1, \text{ if }K\text{ is conjugate to }H\\
0, \text{ otherwise}
\end{cases}$$

Since the dimension of $A(G)$ agrees with the number of conjugacy classes, the collection of $e_H$, with one $H$ chosen from each conjugacy class, forms a basis of $A(G)$.
\medskip

In general, for any $H\subset L\subset G$, define a linear map
$$\varphi_H^L:A(L)\rightarrow\Z$$
which sends each $L$-set to the size of its $H$-fixed subset. We can get a similar basis $\{e_H^L\}$, with one $H$ chosen from each conjugacy class of subgroups of $L$, such that
$$\varphi_K^L(e_H^L)=
\begin{cases}
1, \text{ if }K\text{ is conjugate to }H\text{ in }L\\
0, \text{ otherwise}
\end{cases}$$
\medskip

Now we discuss how these maps interact with transfer and restriction maps:
\smallskip

For any $K\subset L_1\subset L_2\subset G$, the restriction map $R_{L_1}^{L_2}$ keeps each $L_2$-set but views it as an $L_1$-set, hence does not change its $K$-fixed subset. So we have
$$\varphi_K^{L_2}=\varphi_K^{L_1}\circ R_{L_2}^{L_1}.$$

On the other hand, for any $K\subset L_1$, the elements in $\{L_1/K\}$ have isotropy groups conjugate to $K$ inside $L_1$. The elements in $T_{L_1}^{L_2}(\{L_1/K\})=\{L_2/K\}$ have isotropy groups conjugate to $K$ inside $L_2$. There may be more isotropy groups. But these additional isotropy groups are still chosen from the conjugacy class of $K$ inside $G$. According to the construction of the basis element $e_H^L$, we have
$$\varphi_H^{L_2}(T_{L_1}^{L_2}e_K^{L_1})\neq 0\text{ only if }H\text{ is conjugate to a subgroup of }K\text{ in }G$$

Let
$$N(L):=\bigcap_{H\in\mathscr{F},H\subset L}\ker\varphi_H^L$$
In other words, $N(L)$ is generated by $e_H^L$ for all $H\subset L$ and $H\notin\mathscr{F}$. The above discussion tells us that $N(L)$ is closed under transfer and restriction maps. Thus we get a sub-Mackey functor whose value at $G/L$ agrees with $N(L)$. Denote that as $N_\mathscr{F}$.
\smallskip

Notice that for each $L\subset G$, $e_H^L\in N(L)=N_\mathscr{F}(G/L)$ if $H\notin\mathscr{F}$, $e_H^L\in M_\mathscr{F}(G/L)$ if $H\in\mathscr{F}$, according to \autoref{degree0}. Moreover, consider any nontrivial element
$$a_1\{L/K_1\}+a_2\{L/K_2\}+...+a_n\{L/K_n\}\in M_\mathscr{F}(G/L)$$
with $a_1,a_2,...,a_n\neq 0$, $K_1,K_2,...,K_n$ in different conjugacy classes in $\mathscr{F}$, and $|K_1|\leq |K_2|\leq...\leq|K_n|$. The map $\varphi_{K_n}^L$ sends $\{L/K_1\},...,\{L/K_{n-1}\}$ to zero but $\{L/K_n\}$ to a positive value. Thus this element is not in $\ker\varphi_{K_n}^L$, and hence $M_\mathscr{F}(G/L)\cap N(L)=\emptyset$.
\medskip

In conclusion, we have $A_G=M_\mathscr{F}\oplus N_\mathscr{F}$ and $M_\mathscr{F}$ appears as a direct summand. $\Box$

\begin{remark} \label{injectivevar}
When the multiplicative structure is added into consideration, the maps $\varphi_H^L$ appear as the components of the ring isomorphism from $A(L)$ to several copies of $\Z[|G|^{-1}]$, with one $H$ chosen from each conjugacy class inside $L$. The generators $e_H^L$ become idempotent elements. So we can view $M_\mathscr{F}$ as a direct summand of $A_G$ as Green functors. This also explains why $e_HX$ in \autoref{idemsplit} can be expressed as a product with universal spaces. More details can be found in \textbf{[Bar08]} and \textbf{[LMS86}, Chapter V].
\end{remark}

For positive degrees, we have

\begin{theorem} \label{torsion}

For any family $\mathscr{F}$, the $\Z$-graded, Mackey functor valued homology
$$\underline{H}_*^G(E\mathscr{F};A_G)$$
contains only torsion when the degree is positive. Moreover, the torsion only has prime factors which divide $|G|$.

\end{theorem}

This is the most important property we want about $E\mathscr{F}$. We will use the rest of this section to prove it.
\bigskip

The main idea is an induction on the size of $\mathscr{F}$. For the base case, we have

\begin{lemma} \label{orduni}

Let $BG=EG/G$ be the classifying space of principal $G$-bundles. Then $H_*(BG;\Z)$ contains only torsion when the degree is positive. Moreover, the torsion only has prime factors which divide $|G|$.

\end{lemma}

This is a standard result about $EG$. We give one possible proof below. The same idea will be generalized to other universal spaces $E\mathscr{F}$ in \autoref{eqorduni} and \autoref{genuniorb}.

\paragraph{Proof:} Give $EG$ the standard $G$-CW structure which only contains $G$-free cells. Consider the map between cellular chain complexes induced by the projection $EG\rightarrow EG/G=BG$:
$$C_*(EG)\rightarrow C_*(BG).$$
Define another map between chain complexes
$$C_*(BG)\rightarrow C_*(EG)$$
which sends each cell $e$ in $BG$ to the sum of all cells in $EG$ which are sent to $e$ under $EG\rightarrow BG$. Then the composition
$$C_*(BG)\rightarrow C_*(EG)\rightarrow C_*(BG)$$
is the multiplication by $|G|$ since all cells in $EG$ are $G$-free. The same $|G|$-multiplication passes to homology:
$$H_*(BG)\rightarrow H_*(EG)\rightarrow H_*(BG).$$

Since $EG$ is contractible, $H_*(EG)=\Z$, which is concentrated in degree $0$. Thus when the degree is positive, $H_*(BG)$ contains only torsion which divides $|G|$. $\Box$
\bigskip

In the case when $\mathscr{F}$ only contains the trivial subgroup, $E\mathscr{F}$ becomes $EG$ and only has orbit spaces homotopic to $BH$ for $H\subset G$. \autoref{orduni} and \autoref{burnsidehom} imply that it satisfies \autoref{torsion}.
\medskip

For general $\mathscr{F}$, we first give a construction of $E\mathscr{F}$:
\smallskip

\begin{definition} \label{joindef}
The \textbf{join} $X\ast Y$ of two spaces $X,Y$ is defined as
$$X\ast Y:= (X\times Y\times[0,1]\sqcup X\sqcup Y)/\sim,$$
with equivalence relation given by projections
$$X\times Y\times\{0\}\rightarrow X, \text{ }X\times Y\times\{1\}\rightarrow Y.$$
\end{definition}

It's not hard to check:

\begin{lemma} \label{join}
For $G$-CW complexes $X,Y$, $X\ast Y$ has a natural $G$-CW structure. For any subgroup $H$, $(X\ast Y)^H$ is

\textbf{(a)} contractible if and only if at least one of $X^H, Y^H$ is contractible.

\textbf{(b)} empty if and only if both $X^H, Y^H$ are empty.
\end{lemma}

Now in order to construct $E\mathscr{F}$, it suffices to find $G$-spaces $X_1,X_2,...,X_n$, such that:

\textbf{(a)} If $H\in\mathscr{F}$, $(X_i)^H$ is contractible for some $i$.

\textbf{(b)} If $H\notin\mathscr{F}$, $(X_i)^H$ is empty for all $i$.

Then $X_1\ast X_2\ast...\ast X_n$ will be a construction for $E\mathscr{F}$ by the above lemma.
\medskip

For any subgroup $H$, consider the $G$-space $G\times_{N_GH}EW_GH$. Here $N_GH$ is the normalizer of $H$ in $G$. We view $EW_GH$ as a $N_GH$-space. Define the product with $N_GH$ left-acting on $EW_GH$ and right-acting on $G$.

The fixed point subspaces of $G\times_{N_GH}EW_GH$ have the desired properties: For any $a,g\in G$ and $x\in EW_GH$, if the action of $g$ fixes the point $a\times x$:
$$g(a\times x)=ga\times x=a\times x,$$
then there is $h\in N_GH$, such that $hx=x$, $gah^{-1}=a$. Since all points in $EW_GH$ have isotropy group $H$, we have $h\in H$. Then $g=aha^{-1}\in aHa^{-1}$.

Therefore, the isotropy group of $a\times x$ is $aHa^{-1}$, and we have

\begin{lemma} \label{joinpiece}
The fixed point subspace $(G\times_{N_GH}EW_GH)^K$ is non-empty if and only if $K$ is sub-conjugate to $H$. Moreover, when $K=gHg^{-1}$, we have
$$(G\times_{N_GH}EW_GH)^K=\{g\}\times EW_GH,$$
which is contractible.
\end{lemma}

Therefore, $E\mathscr{F}$ can be constructed as the join of $G\times_{N_GH}EW_GH$, with one $H$ chosen from each conjugacy class in $\mathscr{F}$.

\begin{remark} \label{isotropytypes}
Notice that all points in $G\times_{N_GH}EW_GH$ have isotropy groups conjugate to $H$. The collection of all possible isotropy groups of points in $E\mathscr{F}$ agrees with $\mathscr{F}$.
\end{remark}

\medskip

The homology of a join can be computed by the Mayer-Vietoris sequence:
$$X\ast Y=(X\times Y\times[0,1]\sqcup X\sqcup Y)/\sim=(X\times Y\times\left[0,\frac{1}{2}\right]\sqcup X/\sim)\cup(X\times Y\times\left[\frac{1}{2},1\right]\sqcup Y/\sim).$$

We write $X\ast Y$ as the union of two mapping cylinders, which are homotopic to $X,Y$ respectively. The intersection of these two cylinders is $X\times Y\times\{1/2\}\simeq X\times Y$.

Thus we have a short exact sequence
$$0\rightarrow C_*(X\times Y)\rightarrow C_*(X)\oplus C_*(Y)\rightarrow C_*(X\ast Y)\rightarrow 0$$
which passes to fixed point subspaces. So we get a long exact sequence on equivariant homology with coefficients in $A_G$:
$$...\rightarrow \underline{H}_{n+1}^G(X\ast Y;A_G)\rightarrow \underline{H}_n^G(X\times Y;A_G)\rightarrow \underline{H}_n^G(X;A_G)\oplus \underline{H}_n^G(Y;A_G)\rightarrow \underline{H}_n^G(X\ast Y;A_G)\rightarrow...$$

On degree $0$, the map
$$\underline{H}^G_0(X\times Y;A_G)\rightarrow \underline{H}^G_0(X;A_G)\oplus \underline{H}^G_0(Y;A_G)$$
is always an inclusion since the 0th homology is determined by the number of connected components for each fixed point subspace. So we have

\begin{lemma} \label{joinhom}
If $X,Y,X\times Y$ satisfy \autoref{torsion}, then so does $X\ast Y$.
\end{lemma}
\medskip

\paragraph{Proof of \autoref{torsion}:} Assume that \autoref{torsion} is true for all families smaller than $\mathscr{F}$. Let $\mathscr{F}^\prime$ be a smaller family which is obtained by removing the conjugacy class of one largest subgroup $H\in\mathscr{F}$. Then $E\mathscr{F}$ can be constructed as
$$E\mathscr{F}=E\mathscr{F}^\prime\ast(G\times_{N_GH}EW_GH).$$
\medskip

For any $K\lhd L\subset G$, the discussion before \autoref{joinpiece} shows that the space $(G\times_{N_GH}EW_GH)^K$ is either empty or several copies of $EW_GH$. Since all points in each copy of $EW_GH$ have the same isotropy group, $(G\times_{N_GH}W_GH)^K/L$ is either empty or the disjoint union of $BN$ for some subgroups $N\subset N_GH\subset G$. \autoref{orduni} shows that its positive degree homology only contains torsion dividing $|G|$. Thus \autoref{burnsidehom} shows that $G\times_{N_GH}EW_GH$ satisfies \autoref{torsion}.
\medskip

Since $E\mathscr{F}^\prime$ satisfies \autoref{torsion}, according to \autoref{joinhom}, now it suffices to prove the theorem for $E\mathscr{F}^\prime\times(G\times_{N_GH}EW_GH)$. Using $X(H)$ to denote $G\times_{N_GH}EW_GH$, the homology of $E\mathscr{F}^\prime\times X(H)$ can be computed by the equivariant Künneth spectral sequence \textbf{[LM04]}:
$$E_{p,q}^2=\underline{Tor}^{A_G}_{p,q}(\underline{H}_*^G(E\mathscr{F}^\prime;A_G),\underline{H}_*^G(X(H);A_G))\Rightarrow\underline{H}_*^G(E\mathscr{F}^\prime\times X(H);A_G).$$

When we invert $|G|$, the homology of both $E\mathscr{F}^\prime$ and $X(H)$ is concentrated in degree 0. Moreover, $\underline{H}_0^G(E\mathscr{F}^\prime;A_G)=M_{\mathscr{F}^\prime}$ appears as a direct summand of $A_G$ by \autoref{homsplit}. Thus the $E_2$-page collapses into a single box product:
$$E^2=E^2_{0,0}=M_\mathscr{F}\text{ }\Box\text{ } \underline{H}_*^G(X(H);A_G).$$

Therefore, $E\mathscr{F}^\prime\times X(H)$ has trivial homology in positive degrees when $|G|$ is inverted. $\Box$
\bigskip

For special choices of $\mathscr{F}$, the torsion in $\underline{H}_*(E\mathscr{F};A_G)$ may not cover all prime factors of $|G|$. We emphasize the criterion which follows \autoref{burnsidehom}:

\begin{proposition} \label{torsiontype}
The torsion types in $\underline{H}_*^G(E\mathscr{F},A_G)$ come from $H_*(E\mathscr{F}^K/L)$ for all $K\lhd L\subset G$.
\end{proposition}

\section{Splittings of rational $G$-spectra}

In this section, we will reprove \autoref{computesplit} and prove \autoref{orderinvert}. The main idea is to consider the partial splitting induced by a cofiber sequence involving the universal space $E\mathscr{F}$:
$$E\mathscr{F}_+\rightarrow S^0\rightarrow\widetilde{E\mathscr{F}}$$
where the map $E\mathscr{F}_+\rightarrow S^0$ sends $E\mathscr{F}$ into the non-basepoint of $S^0$. The cofiber $\widetilde{E\mathscr{F}}$ can be viewed as the non-equivariant unbased suspension of $E\mathscr{F}$ but with a natural $G$-action.
\smallskip

\begin{theorem} \label{ringsplit}
When certain prime factors of $|G|$ are inverted, this cofiber sequence splits the category of $G$-spectra, in the sense that
$$X\simeq (E\mathscr{F}_+\wedge X)\vee(\widetilde{E\mathscr{F}}\wedge X),$$
$$[E\mathscr{F}_+\wedge X,\widetilde{E\mathscr{F}}\wedge Y]^G=[\widetilde{E\mathscr{F}}\wedge Y,E\mathscr{F}_+\wedge X]^G=0$$
for any $G$-spectra $X,Y$. Moreover, the functors
$$X\mapsto E\mathscr{F}_+\wedge X, \text{ }X\mapsto\widetilde{E\mathscr{F}}\wedge X$$
are symmetric monoidal.
\end{theorem}

We will give a criterion about which prime factors we need to invert in \autoref{invertprime}.
\medskip

The complete splitting in \autoref{computesplit} can be obtained by iterating \autoref{ringsplit} for all choices of $\mathscr{F}$. Instead of full rationalization, we will show that it suffices to invert all prime factors of $|G|$ in our proof.
\bigskip

As based $G$-spaces, both $E\mathscr{F}_+$ and $\widetilde{E\mathscr{F}}$ can be characterized by their fixed point subspaces:
$$E\mathscr{F}_+^H\simeq S^0 \text{ if } H\in\mathscr{F}, \text{ } E\mathscr{F}_+^H\simeq * \text{ if } H\notin\mathscr{F},$$
$$(\widetilde{E\mathscr{F}})^H\simeq S^0 \text{ if } H\notin\mathscr{F}, \text{ } (\widetilde{E\mathscr{F}})^H\simeq * \text{ if } H\in\mathscr{F}.$$

\autoref{monoid} also implies a multiplicative structure on $\widetilde{E\mathscr{F}}$. But unlike $E\mathscr{F}$, the map $S^0\rightarrow\widetilde{E\mathscr{F}}$ also provides a unit, which makes $\widetilde{E\mathscr{F}}$ into a based $G$-Hopf space.
\medskip

Consider the suspension of the cofiber sequence about $E\mathscr{F}$:
$$\Sigma^\infty E\mathscr{F}_+\rightarrow S\rightarrow\Sigma^\infty\widetilde{E\mathscr{F}}.$$

\begin{theorem} \label{cofibsplit}
When $|G|$ is inverted, there exists a left inverse of $\Sigma^\infty E\mathscr{F}_+\rightarrow S$.
\end{theorem}

\paragraph{Proof:} According to \autoref{degree0}, the map $\Sigma^\infty E\mathscr{F}_+\rightarrow S$ induces an inclusion in homology:
$$\underline{HA_G}_0(\Sigma^\infty E\mathscr{F}_+)=:M_\mathscr{F}\hookrightarrow\underline{HA_G}_0S=A_G.$$

When $|G|$ is inverted, \autoref{torsion} shows that the $HA_G$-homology is concentrated in degree 0. Thus it suffices to find a map $S\rightarrow\Sigma^\infty E\mathscr{F}_+$ which induces a projection from $A_G$ to $M_\mathscr{F}$ on their $HA_G$-homology. After composition with $\Sigma^\infty E\mathscr{F}_+\rightarrow S$, we get a self-map of $\Sigma^\infty E\mathscr{F}_+$ that induces an isomorphism on $HA_G$-homology. Since $\Sigma^\infty E\mathscr{F}_+$ is a connective spectrum, this self map must be a weak equivalence.
\smallskip

For each $H\in\mathscr{F}$, the transfer and restriction maps in Mackey functors between images of $G/G$ and $G/H$ are induced by stable maps

$$T: S=\Sigma^\infty S^0=\Sigma^\infty  G/G_+\rightarrow\Sigma^\infty G/H_+$$
$$R: \Sigma^\infty G/H_+\rightarrow \Sigma^\infty G/G_+=S.$$

$R$ is induced by the space-level map $G/H_+\rightarrow S^0$, which sends $G/H$ into the non-basepoint of $S^0$. According to \autoref{isotropytypes}, $G/H$ can be mapped into $E\mathscr{F}$. Thus $R$ is factorized as
$$\Sigma^\infty G/H_+\rightarrow\Sigma^\infty E\mathscr{F}_+\rightarrow S.$$

Let $\iota_H$ be the composition of the first map above and $T$:
$$\iota_H:S\xrightarrow{T}\Sigma^\infty G/H_+\rightarrow\Sigma^\infty E\mathscr{F}_+.$$

The composition
$$S\xrightarrow{\iota_H}\Sigma^\infty E\mathscr{F}_+\rightarrow S$$
agrees with $R\circ T$, which is the equivariant Euler characteristic of $G/H$. We can fully describe the self-map on
$$\underline{HA_G}_0S\cong A_G\cong\underline{\pi}_0S$$
induced by the Euler characteristic as follows:

For each $L\subset G$, the induced self-map on $A_G(G/L)=A(L)$ is the multiplication by $\{G/H\}$ (which is viewed as an $L$-set). We refer to $\textbf{[LMS86}$, V.1 and V.2] for more details.
\medskip

Now we choose one $H$ from each conjugacy class in $\mathscr{F}$. Assign a number $c_H\in\Z[|G|^{-1}]$ for each such $H$. Consider the map
$$\sum_H c_H\iota_H:S\rightarrow\Sigma^\infty E\mathscr{F}_+.$$
The composition
$$S\xrightarrow{\sum_Hc_H\iota_H}\Sigma^\infty E\mathscr{F}_+\rightarrow S$$
induces a self-map on $A_G(G/L)=A(L)$ as the multiplication by $\sum_H c_H\{G/H\}$.

Recall that $\underline{HA_G}_0(\Sigma^\infty E\mathscr{F}_+)(G/L)=M_\mathscr{F}(G/L)\subset A(L)$ is generated by all $\{L/J\}$ with $J\in\mathscr{F}$. The statement that $\sum_Hc_H\iota_H$ induces the projection $A_G\rightarrow M_\mathscr{F}$ is equivalent to
$$\{L/J\}\cdot\sum_Hc_H\{G/H\}=\{L/J\}$$
for any $J\in\mathscr{F}$, $J\subset L$.

According to \autoref{injectivevar}, $\prod_{K\subset L}\varphi_K^L$ is an injective ring homomorphism from $A(L)$ to copies of $\Z$. Thus it suffices to check the equation above on the images under each $\varphi_K^L$. Recall that $\varphi_K^L$ is defined by sending each $L$-set to the number of $K$-fixed points. Thus the equation above under $\varphi_K^L$ becomes
$$|(L/J)^K|\cdot\sum_Hc_H|(G/H)^K|=|(L/J)^K|.$$

Notice that $|(L/J)^K|=0$ for any $K\notin\mathscr{F}$ (since $K$ is not sub-conjugate to $J\in\mathscr{F}$). Thus it suffices to choose $c_H$ such that
$$\sum_Hc_H|(G/H)^K|=\sum_H c_Hs_{(K,H)}=1, \text{ } \forall K\in\mathscr{F}$$
with one $H$ chosen from each conjugacy class in $\mathscr{F}$.

Recall that $s_{(K,H)}$ is non-zero if and only if $K$ is sub-congujate to $H$. Denote this relation by $[K]\leq[H]$. We have

$$\sum_{[K]\leq[H]\subset\mathscr{F}} c_Hs_{(K,H)}=1.$$
Thus
$$c_K=s_{(K,K)}^{-1}\left(1-\sum_{[K]\lneqq[H]\subset\mathscr{F}}c_Hs_{(K,H)}\right).$$

Since $s_{(K,K)}=|W_GK|$ is invertible when $|G|$ is inverted, $c_K$ can be chosen inductively from larger subgroups to smaller ones. Moreover, it's clear that the denominator of any $c_K$ only contains prime factors dividing $|G|$.

Now the map $\sum_Hc_H\iota_H$ gives us the required projection on the 0th homology, hence becomes the left inverse of $\Sigma^\infty E\mathscr{F}_+\rightarrow S$. $\Box$
\bigskip

\paragraph{Proof of \autoref{ringsplit}:} Assume that $|G|$ is inverted. The left inverse in \autoref{cofibsplit} and the multiplicative structure on space $E\mathscr{F}$ make $\Sigma^\infty E\mathscr{F}_+$ into a ring spectrum. So we have a splitting cofiber sequence of ring spectra:
$$\Sigma^\infty E\mathscr{F}_+\rightarrow S\rightarrow\Sigma^\infty\widetilde{E\mathscr{F}}.$$

Thus we have
$$X\simeq(E\mathscr{F}_+\wedge X)\vee(\widetilde{E\mathscr{F}}\wedge X).$$
Since
$$E\mathscr{F}_+\wedge E\mathscr{F}_+\simeq E\mathscr{F}_+\text{ and }\widetilde{E\mathscr{F}}\wedge\widetilde{E\mathscr{F}}\simeq\widetilde{E\mathscr{F}},$$
the functors
$$X\mapsto E\mathscr{F}_+\wedge X,\text{ }X\mapsto\widetilde{E\mathscr{F}}\wedge X$$
are symmetric monoidal.

Moreover, $E\mathscr{F}_+\wedge\widetilde{E\mathscr{F}}\simeq *$ since all its fixed point subspaces are contractible. Use $A,B$ to denote the suspensions of $E\mathscr{F}_+$ and $\widetilde{E\mathscr{F}}$ (in either order). The universal coefficient spectral sequence \textbf{[LM04]} tells us that for any $G$-spectra $X,Y$, we have
$$\underline{Ext}^{*,*}_{\underline{B}_*}(\underline{B}_*(A\wedge X),\underline{(B\wedge Y)}_*)\Rightarrow\underline{(B\wedge Y)}^*(A\wedge X)=[A\wedge X,\Sigma^* B\wedge Y]^G.$$
Since
$$\underline{B}_*(A\wedge X)=\underline{\pi}_*(B\wedge A\wedge X)=0,$$
the spectral sequence has trivial $E_2$-page. Thus we have

$$[E\mathscr{F}_+\wedge X,\widetilde{E\mathscr{F}}\wedge Y]^G=[\widetilde{E\mathscr{F}}\wedge Y,E\mathscr{F}_+\wedge X]^G=0.$$
$\Box$
\medskip

For special choices of $\mathscr{F}$, here is a criterion about which prime factors of $|G|$ we need to invert:

\begin{remark} \label{invertprime}

It suffices to invert the prime factors of the denominators of all $c_H$, and the torsion in $\underline{H}_*^G(E\mathscr{F};A_G)$.

According to \autoref{torsiontype}, the torsion in $\underline{H}_*^G(E\mathscr{F};A_G)$ comes from the torsion in $H_*(E\mathscr{F}^K/L)$ for all $K\lhd L\subset G$.

\end{remark}
\bigskip

Now we list all different families $\mathscr{F}_1,\mathscr{F}_2,...,\mathscr{F}_n$ of $G$. Use \autoref{ringsplit} repeatedly:
$$X\simeq((E\mathscr{F}_1)_+\wedge X)\vee(\widetilde{E\mathscr{F}_1}\wedge X)$$
$$\simeq((E\mathscr{F}_1)_+\wedge(E\mathscr{F}_2)_+\wedge X)\vee((E\mathscr{F}_1)_+\wedge\widetilde{E\mathscr{F}_2}\wedge X)$$
$$\vee(\widetilde{E\mathscr{F}_1}\wedge(E\mathscr{F}_2)_+\wedge X)\vee(\widetilde{E\mathscr{F}_1}\wedge\widetilde{E\mathscr{F}_2}\wedge X)\simeq...$$
Finally we have
$$X\simeq\bigvee_\lambda E\lambda\wedge X$$
where the wedge sum is computed over all functions 
$$\lambda:\{1,2,...,n\}\rightarrow \text{\{based G-spaces\}}$$
sending each $i$ to either $(E\mathscr{F}_i)_+$ or $\widetilde{E\mathscr{F}_i}$. We define $E\lambda$ as the smash product of all elements in the image of $\lambda$.
\medskip

Since the fixed point subspaces of each $E\mathscr{F}_+$ and $\widetilde{E\mathscr{F}}$ are homotopic to either $S^0$ or $*$, the same thing is true for $E\lambda$.

\begin{lemma}
For any $\lambda$, there is at most one conjugacy class $[H]$ such that $(E\lambda)^H\simeq S^0$. Moreover, for each $[H]$, there is a unique $\lambda$ with $(E\lambda)^H\simeq S^0$. We call it $\lambda_H$.
\end{lemma}

\paragraph{Proof:} For any two subgroups $K,H$ which are not in the same conjugacy class, assume that $|K|\leq |H|$. Let $\mathscr{F}$ be the class of all subgroups which are sub-conjugate to $K$. Then
$$(E\mathscr{F}_+)^K\simeq(\widetilde{E\mathscr{F}})^H\simeq S^0,$$
$$(E\mathscr{F}_+)^H\simeq(\widetilde{E\mathscr{F}})^K\simeq *.$$

No matter which of $E\mathscr{F}_+$, $\widetilde{E\mathscr{F}}$ is chosen in $E\lambda$, we have to make one of $(E\lambda)^K,(E\lambda)^H$ contractible.
\smallskip

On the other hand, if $(E\lambda)^H\simeq S^0$, we have to choose $E\mathscr{F}_+$ if $H\in\mathscr{F}$, and $\widetilde{E\mathscr{F}}$ if $H\notin\mathscr{F}$. Thus $\lambda$ is unique. $\Box$
\bigskip

We can view $E\lambda_H$ as another kind of universal space and characterize $E\lambda_H$ by fixed point subspaces:

\begin{proposition} \label{moreclass}
The $K$-fixed subspace of $E\lambda_H$ is homotopic to $S^0$ if $K$ is conjugate to $H$, and contractible otherwise.

Any based $G$-CW complex with the same fixed point spaces is homotopic to $E\lambda_H$.
\end{proposition}

Now we get the splitting in \autoref{computesplit} with just inverting $|G|$. Thus \autoref{orderinvert} is proved.

\section{Splittings for dihedral groups}

The partial splitting in \autoref{ringsplit}, together with \autoref{invertprime}, is quite powerful in the study of $G$-spectra. When we just invert some of the prime factors of $|G|$, the splitting in \autoref{ringsplit} may still happen for special choices of $\mathscr{F}$. We can consider different splittings with different primes inverted in order to recover information about given spectra.

The following three sections are devoted to explicit computations for the dihedral group $G=D_{2p}$ with $p$ an odd prime. We will first compute the $RO(G)$-graded homotopy group $\pi_\bigstar(H\underline{\Z})$, which expresses the equivariant cohomology of a point with coefficients in $\underline{\Z}$. We will improve the computations in \textbf{[KL20]} and \textbf{[Zou18]} and compute $\pi_\bigstar(H\underline{\Z})$ as an $RO(G)$-graded ring with the idea explained above.
\bigskip

Use $\zeta, \tau$ to express the generators in $D_{2p}$ such that $\zeta^p=\tau^2=1$, $\zeta\tau=\tau\zeta^{-1}$.

The $RO(G)$ grading is generated by three different kinds of irreducible representations: The constant representation (denoted by $1$), the sign representation (denoted by $\sigma$), and $p-1$ different 2-dimensional dihedral representations (where $\tau$ acts as a reflection and $\zeta$ acts as a rotation).

The Periodicity Theorem in \textbf{[KL20]} shows that the cellular structures of the representation spheres for different dihedral representations are isomorphic. Let $\gamma$ be any single dihedral representation. It suffices to consider all degrees with form $k+m\sigma+n\gamma$, $\forall k,m,n\in\Z$. The additive and multiplicative structures of homotopy groups on these degrees already reveal all information about $\pi_\bigstar$.
\bigskip

\paragraph{Main Strategy:} Consider the families 
$$\mathscr{F}_1=\{\{1\},\{1,\zeta^i\tau\}:i=1,2,...,p\},$$
$$\mathscr{F}_2=\{\{1\},\{1,\zeta,\zeta^2,...,\zeta^{p-1}\}\}.$$
We will prove:

\begin{proposition} \label{twosplit}
The splitting in \autoref{ringsplit} happens when

\textbf{(a)} $\mathscr{F}=\mathscr{F}_1$ and the prime $p$ is inverted, or

\textbf{(b)} $\mathscr{F}=\mathscr{F}_2$ and the prime $2$ is inverted.
\end{proposition}

Let $A$ be any $G$-ring spectrum. If we can compute the $RO(G)$-graded homotopy and corresponding ring structures, with either $2$ or $p$ inverted, for
$$(E\mathscr{F}_1)_+\wedge A,\text{ }\widetilde{E\mathscr{F}_1}\wedge A, \text{ }(E\mathscr{F}_2)_+\wedge A\text{ and }\widetilde{E\mathscr{F}_2}\wedge A,$$
then we can recover $\pi_\bigstar(A)\otimes\Z[1/2]$ and $\pi_\bigstar(A)\otimes\Z[1/p]$, hence also $\pi_\bigstar(A)$ itself.

Intuitively, we decompose the computation into $C_2$-spectra when $p$ is inverted, and into $C_p$-spectra when $2$ is inverted. Since $\pi_\bigstar(H\underline{Z})$ is already fully computed when $G=C_2$ or $C_p$, we can use these computations as building blocks to recover $\pi_\bigstar(H\underline{\Z})$ for $G=D_{2p}$.
\bigskip

In order to prove \autoref{twosplit}, we need an equivariant generalization of \autoref{orduni}:

\begin{lemma} \label{eqorduni}
Let $X$ be a $D_{2p}$-CW complex which only has $C_p$-free cells. Consider the projection $X\rightarrow X/C_p$, which is a $C_2$-map for any choice of $C_2\subset D_{2p}$ acting on $X$.

\textbf{(a)} The cokernel of 
$$H_*(X;\Z)\rightarrow H_*(X/C_p;\Z)$$
only contains $p$-torsion;

\textbf{(b)} We have an isomorphism
$$H_*(X^{C_2};\Z)\xrightarrow{\cong}H_*((X/C_p)^{C_2};\Z).$$

If $p$ is inverted and the map in part (a) is an injection, then $X\rightarrow X/C_p$ is a $C_2$-equivalence. 
\end{lemma}

\paragraph{Proof:} Part (a) follows by the same argument as in \autoref{orduni}.

For part (b), any $D_{2p}$-cell in $X$ has the form of either $(D_{2p})_+\wedge e$ or $(D_{2p}/C_2)_+\wedge e$. In the first case, its image in $X/C_p$ is not $C_2$-fixed. In the second case, its image is $C_2$-fixed. But no matter which $C_2\subset D_{2p}$ is chosen, the action of $C_2$ fixes exactly one copy of $e$ and exchanges the rest in pairs. Thus we have an isomorphism $C_*(X^{C_2})\cong C_*((X/C_p)^{C_2})$, hence an isomorphism on homology. 

If $p$ is inverted, the cokernel of the map in part (a) disappeared. If that map is also an injection, it becomes an isomorphism. Thus the map $X\rightarrow X/C_p$ induces equivalences on both the underlying spaces and $C_2$-fixed subspaces, hence becomes a $C_2$-equivalence. $\Box$
\medskip

\begin{remark} \label{eqordunimore}
\autoref{eqorduni} also works for $D_{2p}$-CW spectra when we replace the fixed point subspaces by geometric fixed point spectra.

Moreover, since the map $X\wedge Y\rightarrow (X/C_p)\wedge (Y/C_p)$ can be factorized into
$$X\wedge Y\rightarrow (X\wedge Y)/C_p\rightarrow (X/C_p)\wedge(Y/C_p)$$
as $C_2$-maps, the $C_2$-equivalence in \autoref{eqorduni} also preserves multiplicative structures.
\end{remark}
\medskip

\paragraph{Proof of \autoref{twosplit}:} According to \autoref{invertprime} and \autoref{torsiontype}, it suffices to show that the denominators of all $c_H$ and all torsion in $H_*(E\mathscr{F}^K/L)$ (with any $K\lhd L\subset G$) are inverted.
\medskip

As shown in the proof of \autoref{cofibsplit}, $c_H$ is computed by
$$\sum_H c_Hs_{(K,H)}=1,\text{ }\forall K\in\mathscr{F}$$
with one $H$ chosen from each conjugacy class in $\mathscr{F}$. Recall that $s_{(K,H)}=|(G/H)^K|$, which is non-zero if and only if $K$ is sub-conjugate to $H$.

When $\mathscr{F}=\mathscr{F}_1$, we have
$$2pc_{\{1\}}+pc_{\{1,\tau\}}=c_{\{1,\tau\}}=1.$$
Thus $c_{\{1\}}=(1-p)/2p$ and $c_{\{1,\tau\}}=1$. Since $1-p$ is even, the denominators only contain prime $p$.
\smallskip

When $\mathscr{F}=\mathscr{F}_2$, we have
$$2pc_{\{1\}}+2c_{\{1,\zeta,...\}}=2c_{\{1,\zeta,...\}}=1.$$
Thus $c_{\{1\}}=0$ and $c_{\{1,\zeta,...\}}=1/2$, which only have $2$ in denominators.
\bigskip

Up to homotopy, $E\mathscr{F}^K/L$ can be any of $E\mathscr{F}$, $E\mathscr{F}^{C_2}$, $E\mathscr{F}^{C_p}$, $E\mathscr{F}^{D_{2p}}$ $E\mathscr{F}^{C_p}/D_{2p}$, $E\mathscr{F}/C_2$, $E\mathscr{F}/C_p$, $E\mathscr{F}/D_{2p}$.
\bigskip

\textbf{When $\mathscr{F}=\mathscr{F}_1$:} As non-equivariant spaces, we have $E\mathscr{F}\simeq E\mathscr{F}^{C_2}\simeq *$, $E\mathscr{F}^{C_p}=E\mathscr{F}^{D_{2p}}=\emptyset$. Thus $E\mathscr{F}^{C_p}/D_{2p}=\emptyset$.

As a $C_2$-space, $E\mathscr{F}\simeq *$ since both its fixed point subspaces are contractible. Thus $E\mathscr{F}/C_2\simeq *$.

For $E\mathscr{F}/C_p$ and $E\mathscr{F}/D_{2p}$, we can apply \autoref{eqorduni} with $X=E\mathscr{F}$. Notice that $E\mathscr{F}\simeq *$ as a non-equivariant space. Thus the map in \autoref{eqorduni}(a) is an injection. When $p$ is inverted, we have a $C_2$-equivalence $E\mathscr{F}\simeq E\mathscr{F}/C_p$. So we have $E\mathscr{F}/C_p\simeq E\mathscr{F}\simeq *$ and $E\mathscr{F}/D_{2p}\simeq E\mathscr{F}/C_2\simeq *$.

Therefore, the homology groups of these spaces only have $p$-torsion when the degree is positive.
\bigskip

\textbf{When $\mathscr{F}=\mathscr{F}_2$:} We can give an explicit construction as $E\mathscr{F}\simeq EC_2$, with trivial $C_p$-action and the usual $C_2\simeq D_{2p}/C_p$ action.

Thus $E\mathscr{F}^{D_{2p}}=E\mathscr{F}^{C_2}=\emptyset$, $E\mathscr{F}=E\mathscr{F}^{C_p}=E\mathscr{F}/C_p\simeq *$, and
$$E\mathscr{F}^{C_p}/D_{2p}\simeq E\mathscr{F}/C_2\simeq E\mathscr{F}/D_{2p}\simeq BC_2.$$
According to \autoref{orduni}, the homology groups of all these spaces only have $2$-torsion when the degree is positive.
\bigskip

In conclusion, we have the splitting when $\mathscr{F}=\mathscr{F}_1$ and $p$ inverted, or when $\mathscr{F}=\mathscr{F}_2$ and $2$ inverted. $\Box$
\bigskip

Now it remains to compute the $RO(G)$-graded homotopy for
$$(E\mathscr{F}_1)_+\wedge A,\text{ }\widetilde{E\mathscr{F}_1}\wedge A, \text{ }(E\mathscr{F}_2)_+\wedge A\text{ and }\widetilde{E\mathscr{F}_2}\wedge A$$
with $A=H\underline{\Z}$.
\medskip

\section{Computation of $\pi_\bigstar(H\underline{\Z})$}

We will use the same way to define generators as in \textbf{[HHR17]} Definition 3.4:

\begin{definition} \label{generator}

For any actual $G$-representation $V$ with $V^G=0$, let $a_V\in\pi_{-V}(S^0)$ be the map $S^0\rightarrow S^V$ embedding $S^0$ to $0$ and $\infty$. We also use $a_V$ to denote its Hurewicz image in $\pi_{-V}(H\underline{\Z})$.

For any actual orientable representation $V$ of dimension $n$, let $u_V$ be the generator of $\pi_{n-V}(H\underline{\Z})=H_n^G(S^V;\underline{\Z})$ which restricts to the choice of orientation in
$$\underline{H}_n^G(S^V;\underline{\Z})(G/e)\cong H_n(S^n;\Z).$$

\end{definition}

Some important relations on these generators are given below:

\begin{proposition} \label{genrelation}
\textbf{(a)} For any $V_1,V_2$,
$$a_{V_1+V_2}=a_{V_1}a_{V_2},\text{ }u_{V_1+V_2}=u_{V_1}u_{V_2}.$$

\textbf{(b)} Let $G_V$ be the isotropy subgroup of $V$. Then $|G/G_V|a_V=0$.

\textbf{(c)} For $V,W$ both oriented with dimension $2$, with $G_V\subset G_W$, we have
$$a_Wu_V=|G_W/G_V|a_Vu_W.$$
\end{proposition}
\medskip

We consider the cases $G=C_2$ and $G=C_p$ first:

For $G=C_2$, we have the trivial representation and the sign representation $\sigma$.

For $G=C_p$, we have the trivial representation and $p-1$ different $2$-dimensional rotation representations. Again, the periodicity theorem in $[\textbf{KL20}]$ works since any two such rotation representation spheres have isomorphic cellular structures. It suffices to consider only one $2$-dimensional rotation, which we denote as $\lambda$.
\medskip

The computations of $\pi_\bigstar(H\underline{\Z})$ for these two cases can be found in \textbf{[Zen18]} and are listed below:

\begin{theorem} \label{primecase}
For $G=C_2$,
$$\pi_\bigstar(H\underline{\Z})=\Z[u_{2\sigma},a_\sigma]/(2a_\sigma)\oplus \left(\bigoplus_{i>0}2\Z\langle u_{2\sigma}^{-i}\rangle\right)\oplus\left(\bigoplus_{j,k>0}\Z/2\langle\Sigma^{-1}u_{2\sigma}^{-j}a_\sigma^{-k}\rangle\right).$$

For $G=C_p$,
$$\pi_\bigstar(H\underline{\Z})=\Z[u_\lambda,a_\lambda]/(pa_\lambda)\oplus\left(\bigoplus_{i>0}p\Z\langle u_\lambda^{-i}\rangle\right)\oplus\left(\bigoplus_{j,k>0}\Z/p\langle\Sigma^{-1}u_\lambda^{-j}a_\lambda^{-k}\rangle\right).$$
\end{theorem}
\medskip

Now let $G=D_{2p}$.

\subsection{Computation of $(E\mathscr{F}_1)_+\wedge H\underline{\Z}$}

We have $\mathscr{F}=\mathscr{F}_1$ and $p$ inverted in this case.
\medskip

For any virtual representation $V$, since $E\mathscr{F}_+\wedge S^V$ is $C_p$-free and $H\underline{\Z}$ is split, \autoref{ringsplit} and the Adams isomorphism tell us
$$[S^V,E\mathscr{F}_+\wedge H\underline{\Z}]^G\cong [E\mathscr{F}_+\wedge S^V, H\underline{\Z}]^G\cong [(E\mathscr{F}_+\wedge S^V)/C_p,H\underline{\Z}]^{C_2}.$$
We refer to \textbf{[May96}, XVI.2] for an explanation of split equivariant spectra.

Write $V=k\gamma+V_0$, where $V_0$ has no copies of $\gamma$. Consider the following composition:
$$E\mathscr{F}_+\wedge S^V=E\mathscr{F}_+\wedge(S^{sign(k)\gamma})^{\wedge|k|}\wedge S^{V_0}\rightarrow(E\mathscr{F}_+\wedge(S^{sign(k)\gamma})^{\wedge|k|}\wedge S^{V_0})/C_p=(E\mathscr{F}\wedge S^V)/C_p$$
$$\rightarrow (E\mathscr{F}_+/C_p)\wedge (S^{sign(k)\gamma}/C_p)^{\wedge|k|}\wedge(S^{V_0}/C_p).$$
Here $sign(k)=k/|k|$ if $k\neq 0$, and $sign(0)=0$.

Since $E\mathscr{F}_+\simeq (EC_p)_+$ as a $C_p$-space, when $p$ is inverted, \autoref{orduni} gives us a non-equivariant equivalence $E\mathscr{F}_+\simeq (E\mathscr{F})_+/C_p$.

The map $S^{sign(k)\gamma}\rightarrow S^{sign(k)\gamma}/C_p$ induces multiplication by $p^{sign(k)}$ on non-equivariant homology, hence is also a non-equivariant equivalence.

Moreover, since $C_p$ acts trivially on $S^{V_0}$, we have $S^{V_0}=S^{V_0}/C_p$. So the composition above is a non-equivariant equivalence. We know that
$$H\Z_*(E\mathscr{F}_+\wedge S^V)\rightarrow H\Z_*((E\mathscr{F}_+\wedge S^V)/C_p)$$
is an injection. According to \autoref{eqorduni} and \autoref{eqordunimore}, 
$$E\mathscr{F}_+\wedge S^V\rightarrow(E\mathscr{F}_+\wedge S^V)/C_p$$
is an equivalence of $C_2$-spectra. So we have
$$[S^V,E\mathscr{F}_+\wedge H\underline{\Z}]^G\cong [(E\mathscr{F}_+\wedge S^V)/C_p,H\underline{\Z}]^{C_2}\cong [E\mathscr{F}_+\wedge S^V, H\underline{\Z}]^{C_2}\cong [S^V,H\underline{\Z}]^{C_2}$$
where the last isomorphism is induced by the fact that $E\mathscr{F}\simeq *$ as a $C_2$-space.

To be more precise:

\begin{proposition} \label{z11}
For any integers $k,m,n$,
$$[S^{k+m\sigma+n\gamma},E\mathscr{F}_+\wedge H\underline{\Z}]^G\cong[S^{(k+n)+(m+n)\sigma},H\underline{\Z}]^{C_2}.$$
Moreover, this isomorphism preserves the multiplicative structures.
\end{proposition}

The isomorphism above induces a ring map from the $RO(D_{2p})$-graded ring $\pi_\bigstar^G(E\mathscr{F}_+\wedge H\underline{\Z})$ to the $RO(C_2)$-graded ring $\pi_\bigstar^{C_2}(H\underline{\Z})$. We can compute $\pi_\bigstar^G(E\mathscr{F}_+\wedge H\underline{\Z})$ by tracing the generators in the first equation of \autoref{primecase}. We have:
\medskip

\textbf{(1)}
$$\pi_0^{C_2}(H\underline{\Z})\cong\pi_{c(1+\sigma-\gamma)}^G(E\mathscr{F}_+\wedge H\underline{\Z}),\text{ }\forall c\in\Z.$$
We use $u_{\gamma-\sigma}$ to denote the generator of $\pi_{1+\sigma-\gamma}^G(E\mathscr{F}_+\wedge H\underline{\Z})$. Then $u_{\gamma-\sigma}$ is invertible.
\smallskip

\textbf{(2)}
$$\pi_{k+m\sigma}^{C_2}(H\underline{\Z})\cong\pi_{k+m\sigma+c(1+\sigma-\gamma)}^G(E\mathscr{F}_+\wedge H\underline{\Z}),\text{ }\forall c\in\Z.$$
Thus the pre-image of $u_{2\sigma},a_\sigma\in\pi_\bigstar^{C_2}(H\underline{\Z})$ are $u_{2\sigma}u_{\gamma-\sigma}^c,a_\sigma u_{\gamma-\sigma}^c\in \pi_\bigstar^G(E\mathscr{F}_+\wedge H\underline{\Z})$, for any $c\in\Z$.
\medskip

Therefore, $\pi_\bigstar^G(E\mathscr{F}_+\wedge H\underline{\Z})$ can be obtained by adding the invertible $u_{\gamma-\sigma}$ to $\pi_\bigstar^{C_2}(H\underline{\Z})$:
$$\pi_\bigstar(E\mathscr{F}_+\wedge H\underline{\Z})[\frac{1}{p}]=\Z[\frac{1}{p}][u_{2\sigma},a_\sigma,u_{\gamma-\sigma}^\pm]/(2a_\sigma)$$
$$\oplus\left(\bigoplus_{i>0}2\Z[\frac{1}{p}][u_{\gamma-\sigma}^\pm]\langle u_{2\sigma}^{-i}\rangle\right)\oplus\left(\bigoplus_{j,k>0}\Z/2[u_{\gamma-\sigma}^\pm]\langle\Sigma^{-1}u_{2\sigma}^{-j}a_\sigma^{-k}\rangle\right).$$

\subsection{Computation of $(E\mathscr{F}_2)_+\wedge H\underline{\Z}$}

We have $\mathscr{F}=\mathscr{F}_2$ and $2$ inverted in this case. We use the explicit construction of $E\mathscr{F}$ as $E(D_{2p}/C_p)=EC_2$.

The Adams isomorphism does not work since $C_2$ is not a normal subgroup of $D_{2p}$. Instead, we construct a collapsing spectral sequence.
\medskip

Give $(EC_2)_+$ the standard cellular structure with one $G$-cell $(G/C_p)_+\wedge e_n$ for each $n=0,1,2,...$, with boundary maps
$$e_0\xleftarrow{1-\tau}e_1\xleftarrow{1+\tau}e_2\xleftarrow{1-\tau}...$$

For any virtual representation $V$, give
$$[S^V,E\mathscr{F}_+\wedge H\underline{\Z}]^G_*\cong [E\mathscr{F}_+\wedge S^V,E\mathscr{F}_+\wedge H\underline{\Z}]^G_*\cong [E\mathscr{F}_+\wedge S^V, H\underline{\Z}]^G_*$$
a filtration induced by the skeleton of $E\mathscr{F}$. This filtration implies a cohomological spectral sequence, whose $E_1$-page, together with the boundary map, can be expressed as the following chain complex:
$$[(G/C_p)_+\wedge S^V,H\underline{\Z}]^G_*\xrightarrow{1-\tau}[(G/C_p)_+\wedge S^V,H\underline{\Z}]^G_*\xrightarrow{1+\tau}[(G/C_p)_+\wedge S^V,H\underline{\Z}]^G_*\xrightarrow{1-\tau}....$$
The action of $\tau$ on $[(G/C_p)_+\wedge S^V,H\underline{\Z}]^G_*$ sends each map $f:(G/C_p)_+\wedge S^V\rightarrow H\underline{\Z}$ to the composition
$$(G/C_p)_+\wedge S^V\xrightarrow{\tau\wedge id} (G/C_p)_+\wedge S^V\xrightarrow{f} H\underline{\Z}.$$

Since 
$$[(G/C_p)_+\wedge S^V,H\underline{\Z}]^G_*\cong [S^V,H\underline{\Z}]^{C_p}_*,$$
the sequence above becomes
$$[S^V,H\underline{\Z}]^{C_p}_*\xrightarrow{1-\tau}[S^V,H\underline{\Z}]^{C_p}_*\xrightarrow{1+\tau}[S^V,H\underline{\Z}]^{C_p}_*\xrightarrow{1-\tau}....$$
The action of $\tau$ on $[S^V,H\underline{\Z}]^{C_p}_*$ is given by conjugation.
\medskip

We will show that the action of $\tau$ is multiplication by either $1$ or $-1$ on each $[S^V,\Sigma^tH\underline{\Z}]^{C_p}$. First we consider the case that $V$ has no copies of $\sigma$:
\medskip

\begin{proposition} \label{tauactionp}
Let $V=k+n\gamma$ for some integers $k,n$. Then the action of $\tau$ on $[S^V,\Sigma^t H\underline{\Z}]^{C_p}$ is multiplication by $-1$ if $|k-t|\equiv 2$ or $3$ (mod $4$). Otherwise, the action is multiplication by $1$.
\end{proposition}

\paragraph{Proof:} Notice that $[S^V,\Sigma^t H\underline{\Z}]^{C_p}$ can be viewed as the equivariant cohomology or homology of $S^{|n|\gamma}$ with coefficients in $\underline{\Z}$. The computation can be done in a cellular way. We first give $S^{|n|\gamma}$ a $C_p$-CW structure.
\smallskip

We view $S^\gamma$ as $\gamma$ compactified at $\infty$. The action of $\zeta$ is the counter-clockwise rotation by $2\pi/p$ and the action of $\tau$ is the reflection by the $x$-axis. The standard $CW$ structure of $S^\gamma$ can be described as follows:
\medskip

As a based non-equivariant space, $S^\gamma\cong S^2$ can be constructed by one single $2$-cell and the base point at the origin. Denote that $2$-cell by $a$.
\smallskip

As a based $C_p$-space, $S^\gamma\cong S^\lambda$ can be constructed by:

$2$-cells $b_1,b_2,...,b_p$, such that each $b_i$ has image
$$\left\{(r\cos\theta,r\sin\theta):\text{ }0\leq r\leq\infty, \frac{2\pi(i-1)}{p}\leq\theta\leq\frac{2\pi i}{p}\right\};$$

$1$-cells $c_1,c_2,...,c_p$, such that each $c_i$ as image
$$\left\{\left(r\cos\frac{2\pi i}{p},r\sin\frac{2\pi i}{p}\right):\text{ }0\leq r\leq\infty\right\};$$

One $0$-cell $d$ at $\infty$ and the base point at $0$.

Here the subscripts of $b_i$ and $c_i$ are defined modulo $p$.
\medskip

For both $CW$ structures, the action of $\tau$ on $S^\gamma$ can be made into a cellular map:
\smallskip

The action of $\tau$ on the whole $S^\gamma$ is a reflection. Thus $\tau a=-a$.

The action of $\tau$ gives a permutation among the images of $b_1,b_2,...,b_p$ and reverses the orientations of these $2$-cells. Thus $\tau b_i=-b_{p+1-i}$.

The action of $\tau$ gives a permutation among the images of $c_1,c_2,...,c_p$ and keeps the orientations of these $1$-cells. Thus $\tau c_i=c_{p+1-i}$.

Finally, $\tau d=d$ since the point $\infty$ is fixed.
\medskip

Since $S^{|n|\gamma}=S^\gamma\wedge S^\gamma\wedge...\wedge S^\gamma$, we can construct a CW structure as follows:

For $i\in\Z$ and $0\leq j\leq |n|-1$, let
$$e_{2j+2,i}:=(\zeta^ia,...,\zeta^ia,b_i,d,...,d)$$
$$e_{2j+1,i}:=(\zeta^ia,...,\zeta^ia,c_i,d,...,d)$$
with $j$ copies of $\zeta^ia$ and $|n|-1-j$ copies of $d$.

We also have one $0$-cell $e_0:=(d,...,d)$ and the base point.

Notice that $\zeta b_i=b_{i+1}$, $\zeta c_i=c_{i+1}$. This construction makes $S^{|n|\gamma}$ into a $C_p$-CW complex with one $C_p$-cell (except the base point) in each degree between $0$ and $2|n|$. The action of $\tau$ on $S^{|n|\gamma}$ is induced by the action on each copy of $S^\gamma$, hence is also made cellular:
$$\tau e_{2j+2,i}=(-1)^{j+1}e_{2j+2,p+1-i},\text{ }\tau e_{2j+1,i}=(-1)^je_{2j+1,p+1-i},\text{ }\tau e_0=e_0.$$
\medskip

Recall that $H_*^{C_p}(S^{|n|\gamma};\underline{\Z})$ and $H^*_{C_p}(S^{|n|\gamma};\underline{\Z})$ can be computed by the following chain and cochain complexes:
$$C_*^{C_p}(S^{|n|\gamma};\underline{\Z}):=\underline{C}_*(S^{|n|\gamma})\otimes_{\mathcal{O}_{C_p}}\underline{\Z};$$
$$C^*_{C_p}(S^{|n|\gamma};\underline{\Z}):=Hom_{\mathcal{O}_{C_p}}(\underline{C}_*(S^{|n|\gamma});\underline{\Z}),$$
where we view $\underline{\Z}$ as a covariant coefficient system in the first equation and a contravariant one in the second equation.

The actions of $\tau$ on the homology and cohomology are obtained by applying $\tau$ on both $\underline{C}_*(S^{|n|\gamma})$ and $\underline{\Z}$.
\smallskip

Since $S^{|n|\gamma}$ has only $C_p$-free cells in positive degrees and only $C_p$-fixed cells in degree $0$, and $\underline{\Z}$ is the constant coefficient system, we have
$$C_*^{C_p}(S^{|n|\gamma};\underline{\Z}):=\underline{C}_*(S^{|n|\gamma})\otimes_{\mathcal{O}_{C_p}}\underline{\Z}\cong C_*(S^{|n|\gamma})/C_p;$$
$$C^*_{C_p}(S^{|n|\gamma};\underline{\Z}):=Hom_{\mathcal{O}_{C_p}}(\underline{C}_*(S^{|n|\gamma});\underline{\Z})\cong Hom_\Z(C_*(S^{|n|\gamma})/C_p;\Z).$$

Let $e_{2j}$ and $e_{2j-1}$ be the orbits of cells $e_{2j,i}$ and $e_{2j-1,i}$. The induced $\tau$-action on $C_*(S^{|n|\gamma})/C_p$ is expressed as
$$\tau e_{2j}=(-1)^je_{2j},\text{ }\tau e_{2j-1}=(-1)^{j-1}e_{2j-1}.$$

Therefore, the $\tau$-action on the homology and cohomology is multiplication by $-1$ when the degree $\equiv$ $2$ or $3$ (mod $4$). Otherwise the $\tau$-action is multiplication by $1$. $\Box$
\bigskip

Now we consider the case that $V$ also contains copies of $\sigma$. For any $D_{2p}$-spectra $X,Y$, let $\tau$ act on $[X,Y]^{C_p}$ by conjugation. Since $S^\sigma\cong S^1$ as $C_p$-spaces, we have
$$[X,Y]^{C_p}\cong [X,S^{m(\sigma-1)}\wedge Y]^{C_p}.$$
Since $\tau$ acts on $S^\sigma$ as a reflection, the composition
$$[X,Y]^{C_p}\cong [X,S^{m(\sigma-1)}\wedge Y]^{C_p}\xrightarrow{\tau}[X,S^{m(\sigma-1)}\wedge Y]^{C_p}\cong [X,Y]^{C_p}$$
agrees with the $\tau$-action on $[X,Y]^{C_p}$ multiplied by $(-1)^m$.
\medskip

Together with \autoref{tauactionp}, we get

\begin{proposition} \label{tauaction}
Let $V=k+m\sigma+n\gamma$. The action of $\tau$ on $[S^V,\Sigma^tH\underline{\Z}]^{C_p}$ is multiplication by $-1$ if $[|k+m-t|/2]+m$ is odd. Otherwise, the action is multiplication by $1$.
\end{proposition}
\bigskip

Recall that in the spectral sequence computing $[S^V,E\mathscr{F}_+\wedge H\underline{\Z}]^G_*$, the $E_1$-page and the boundary maps are given by
$$[S^V,\Sigma^tH\underline{\Z}]^{C_p}\xrightarrow{1-\tau}[S^V,\Sigma^tH\underline{\Z}]^{C_p}\xrightarrow{1+\tau}[S^V,\Sigma^tH\underline{\Z}]^{C_p}\xrightarrow{1-\tau}...,$$
which can be written as either
$$[S^V,\Sigma^tH\underline{\Z}]^{C_p}\xrightarrow{0}[S^V,\Sigma^tH\underline{\Z}]^{C_p}\xrightarrow{2}[S^V,\Sigma^tH\underline{\Z}]^{C_p}\xrightarrow{0}....$$
or
$$[S^V,\Sigma^tH\underline{\Z}]^{C_p}\xrightarrow{2}[S^V,\Sigma^tH\underline{\Z}]^{C_p}\xrightarrow{0}[S^V,\Sigma^tH\underline{\Z}]^{C_p}\xrightarrow{2}....$$
depending on whether the $\tau$-action becomes the multiplication by $1$ or $-1$.

Since $2$ is inverted, the $E_2$-page vanishes except in the bottom line. Consider the case of $t=0$. We can compute the homotopy of $E\mathscr{F}_+\wedge H\underline{\Z}$ by this $E_2$-page:

\begin{proposition} \label{z12}
When $2$ is inverted, for any integers $k,m,n$,
$$[S^{k+m\sigma+n\gamma},E\mathscr{F}_+\wedge H\underline{\Z}]^G\cong
\begin{cases}
0,\text{ if }[|k+m|/2]+m\text{ is odd,}\\
[S^{(k+m)+n\lambda},H\underline{\Z}]^{C_p},\text{ if }[|k+m|/2]+m\text{ is even.}
\end{cases}$$
Moreover, this isomorphism preserves the multiplicative structures.
\end{proposition}

The isomorphism in the second case is obtained by sending each $G$-map
$$S^V\rightarrow E\mathscr{F}_+\wedge H\underline{\Z}$$
to the underlying $C_p$-map of
$$S^V\rightarrow E\mathscr{F}_+\wedge H\underline{\Z}\rightarrow S\wedge H\underline{\Z}\simeq H\underline{\Z}.$$
Thus the multiplicative structures are preserved.

\begin{remark}
According to the second equation in \autoref{primecase}, $\pi_{a+b\lambda}^{C_p}(H\underline{\Z})$ is non-zero only if $a$ is positive even, or negative odd. So we can rewrite \autoref{z12} as
$$[S^{k+m\sigma+n\gamma},E\mathscr{F}_+\wedge H\underline{\Z}]^G\cong
\begin{cases}
0,\text{ if }[(k+3m+1)/2]\text{ is odd,}\\
[S^{(k+m)+n\lambda},H\underline{\Z}]^{C_p},\text{ if }[(k+3m+1)/2]+m\text{ is even.}
\end{cases}$$
\end{remark}

The equation above induces a ring map from the $RO(D_{2p})$-graded ring $\pi_\bigstar^G(E\mathscr{F}_+\wedge H\underline{\Z})$ to the $RO(C_p)$-graded ring $\pi_\bigstar^{C_p}(H\underline{\Z})$. Similarly, in order to compute $\pi_\bigstar^G(E\mathscr{F}_+\wedge H\underline{\Z})$, we only need to trace the generators in the second equation in \autoref{primecase}:
\medskip

\textbf{(1)}
$$\pi_0^{C_p}(H\underline{\Z})\cong\pi_{2c(1-\sigma)}^G(E\mathscr{F}_+\wedge H\underline{\Z}),\text{ }\forall c\in\Z.$$
Thus $u_{2\sigma}$, as the generator of $\pi_{1-\sigma}^G(E\mathscr{F}_+\wedge H\underline{\Z})$, becomes invertible.
\smallskip

\textbf{(2)}
$$\pi_{2-\lambda}^{C_p}(H\underline{\Z})\cong\pi_{1+\sigma-\gamma+2c(1-\sigma)}^G(E\mathscr{F}_+\wedge H\underline{\Z}),\text{ }\forall c\in\Z.$$
Thus the pre-image of $u_\lambda\in\pi_\bigstar^{C_p}(H\underline{\Z})$ is $u_{\gamma-\sigma}u_{2\sigma}^c$, for any $c\in\Z$.
\smallskip

\textbf{(3)}
$$\pi_{-\lambda}^{C_p}(H\underline{\Z})\cong\pi_{-\gamma+2c(1-\sigma)}^G(E\mathscr{F}_+\wedge H\underline{\Z}),\text{ }\forall c\in\Z.$$
Thus the pre-image of $a_\lambda\in\pi_\bigstar^{C_p}(H\underline{\Z})$ is $a_\gamma u_{2\sigma}^c$, for any $c\in\Z$.
\medskip

Therefore, we have
$$\pi_\bigstar(E\mathscr{F}_+\wedge H\underline{\Z})[\frac{1}{2}]=\Z[\frac{1}{2}][u_{\gamma-\sigma},a_\gamma,u_{2\sigma}^\pm]/(pa_\gamma)$$
$$\oplus\left(\bigoplus_{i>0}p\Z[\frac{1}{2}][u_{2\sigma}^\pm]\langle u_{\gamma-\sigma}^{-i}\rangle\right)\oplus\left(\bigoplus_{j,k>0}\Z/p[u_{2\sigma}^\pm]\langle\Sigma^{-1}u_{\gamma-\sigma}^{-j}a_\gamma^{-k}\rangle\right).$$

\subsection{Computation of $\widetilde{E\mathscr{F}_1}\wedge H\underline{\Z}$ and $\widetilde{E\mathscr{F}_2}\wedge H\underline{\Z}$}

In these two cases, we can show that
$$\widetilde{E\mathscr{F}_1}\wedge H\underline{\Z}\simeq \widetilde{E\mathscr{F}_2}\wedge H\underline{\Z}\simeq *.$$
\medskip

We can compute $\underline{\pi}_*(\widetilde{E\mathscr{F}}\wedge H\underline{\Z})$ by the universal coefficient spectral sequence in \textbf{[LM04]}:
$$\underline{Tor}_{*,*}^{A_G}(\underline{HA_G}_*\widetilde{E\mathscr{F}},\underline{\Z})\Rightarrow\underline{H\underline{\Z}}_*\widetilde{E\mathscr{F}}=\underline{\pi}_*(\widetilde{E\mathscr{F}}\wedge H\underline{\Z}).$$

When either $2$ or $p$ is inverted and the splitting in \autoref{ringsplit} happens, $\underline{HA_G}_*\widetilde{E\mathscr{F}}$ is concentrated in degree $0$, and becomes the complement Green functor of $\underline{HA_G}_*(E\mathscr{F}_+)$ in $A_G$. We use the same notations $M_\mathscr{F}$ and $N_\mathscr{F}$ as in section 2 to denote $\underline{HA_G}_*(E\mathscr{F}_+)$ and $\underline{HA_G}_*\widetilde{E\mathscr{F}}$.

Recall that, as a sub-Green functor of $A_G$, the value of $N_\mathscr{F}$ on $G/L$ is generated by all $\{L/K\}$ with $K\notin\mathscr{F}$.

Since $N_\mathscr{F}$ is a direct summand of $A_G$,
$$\underline{Tor}_{s,t}^{A_G}(N_\mathscr{F},
\underline{\Z})=
\begin{cases}
\underline{\Z}\text{ }\Box\text{ }N_\mathscr{F},\text{ if }s=t=0,\\
0,\text{ otherwise.}
\end{cases}$$
Thus
$$\underline{\pi}_*(\widetilde{E\mathscr{F}}\wedge H\underline{\Z})=\underline{\Z}\text{ }\Box\text{ }N_\mathscr{F},$$
which is concentrated in degree $0$.

The box product can be computed as a coend, for which we refer to \textbf{[LM04]}. In the case when $G$ is finite, we have
$$(M\text{ }\Box\text{ }N)(G/H)=\left(\bigoplus_{K\subset H}M(G/K)\otimes N(G/K)\right)/\sim,$$
where the equivalence relation is generated by 
$$a\otimes T_{K_1}^{K_2}(b)\sim R_{K_1}^{K_2}(a)\otimes b, \text{ } T_{K_1}^{K_2}(b)\otimes a\sim b\otimes R_{K_1}^{K_2}(a)$$
for any $K_1\subset K_2$, and
$$M(G/K_1)\otimes N(G/K_1)\cong M(G/K_2)\otimes N(G/K_2)$$
for all conjugacy pairs $K_1,K_2$.
\bigskip

When $\mathscr{F}=\mathscr{F}_1$ and $p$ is inverted, $N_\mathscr{F}(G/G)$ is generated by $\{G/G\}$ and $\{G/C_p\}$, $N_\mathscr{F}(G/C_p)$ is generated by $\{C_p/C_p\}$, and we have
$$N_\mathscr{F}(G/e)=N_\mathscr{F}(G/C_2)=0.$$
Thus
$$(\underline{\Z}\text{ }\Box\text{ }N_\mathscr{F})(G/e)=(\underline{\Z}\text{ }\Box\text{ }N_\mathscr{F})(G/C_2)=0.$$

For any $a\in\underline{\Z}(G/e)$ and $b\in N_\mathscr{F}(G/C_p)$, we have
$$Tr_e^{C_p}(a)\otimes b\sim a\otimes R_e^{C_p}(b)=a\otimes 0=0.$$
However, $Tr_e^{C_p}$ in $\underline{\Z}$ is the multiplication by $p$, which is invertible. Thus this equivalence relation identifies the whole $\underline{\Z}(G/C_p)\otimes N_\mathscr{F}(G/C_p)$ to zero.

By a similar argument, the equivalence relation between $G/G$ and $G/C_2$ identifies the whole $\underline{\Z}(G/G)\otimes N_\mathscr{F}(G/G)$ to zero. Thus we have $\underline{\Z}\text{ }\Box\text{ }N_\mathscr{F}=0$.
\medskip

When $\mathscr{F}=\mathscr{F}_2$ and $2$ is inverted, we just need to exchange $2$ and $p$ in the above argument. So we still have $\underline{\Z}\text{ }\Box\text{ }N_\mathscr{F}=0$.
\medskip

Therefore, 
$$\widetilde{E\mathscr{F}_1}\wedge H\underline{\Z}\simeq \widetilde{E\mathscr{F}_2}\wedge H\underline{\Z}\simeq *.$$

\subsection{Final result}

We already computed $\pi_\bigstar(H\underline{\Z})$ with either $2$ or $p$ inverted. Now we just need to glue them together:

\begin{theorem} \label{hz}
For $G=D_{2p}$,
$$\pi_\bigstar(H\underline{\Z})=\Z[u_{\gamma-\sigma},u_{2\sigma},a_\sigma,a_\gamma]/(2a_\sigma,pa_\gamma)$$
$$\oplus\left(\bigoplus_{i>0}2\Z[u_{\gamma-\sigma}]\langle u_{2\sigma}^{-i}\rangle\right)\oplus\left(\bigoplus_{i>0}p\Z[u_{2\sigma}]\langle u_{\gamma-\sigma}^{-i}\rangle\right)\oplus\left(\bigoplus_{j,k>0}2p\Z\langle u_{2\sigma}^{-j}u_{\gamma-\sigma}^{-k}\rangle\right)$$
$$\oplus\left(\bigoplus_{j,k>0}\Z/p[u_{\gamma-\sigma}]\langle a_\gamma^j u_{2\sigma}^{-k}\rangle\right)\oplus\left(\bigoplus_{j,k>0}\Z/2[u_{2\sigma}]\langle a_\sigma^j u_{\gamma-\sigma}^{-k}\rangle\right)$$
$$\oplus\left(\bigoplus_{j,k>0}\Z/2[u_{\gamma-\sigma}^\pm]\langle\Sigma^{-1}u_{2\sigma}^{-j}a_\sigma^{-k}\rangle\right)\oplus\left(\bigoplus_{j,k>0}\Z/p[u_{2\sigma}^\pm]\langle\Sigma^{-1}u_{\gamma-\sigma}^{-j}a_\gamma^{-k}\rangle\right).$$

\end{theorem}

\section{Applications to general $D_{2p}$-spectra}
Let $G=D_{2p}$. In this section, we will generalize the computations in section 5 to other $G$-spectra. For any $G$-spectra $A$, \autoref{twosplit} shows that the information of $\pi_\bigstar(A)$ is covered by the homotopy of
$$(E\mathscr{F}_1)_+\wedge A,\text{ }\widetilde{E\mathscr{F}_1}\wedge A, \text{ }(E\mathscr{F}_2)_+\wedge A\text{ and }\widetilde{E\mathscr{F}_2}\wedge A$$
with either $2$ or $p$ inverted.

The main results of this section are given below:

\begin{theorem} \label{generald2p}
Let $k,m,n$ be arbitrary integers.

\textbf{(a)} Let $A$ be a split $G$-spectra. When $p$ is inverted, we have
$$[S^{k+m\sigma+n\gamma},(E\mathscr{F}_1)_+\wedge A]^G\cong[S^{(k+n)+(m+n)\sigma},A^{C_p}]^{C_2}.$$

\textbf{(b)} Let $A=HM$ for some Mackey functor $M$ such that $M(G/H)$ has a trivial $W_GH$-action for all $H\subset G$. When $2$ is inverted, we have
$$[S^{k+m\sigma+n\gamma},(E\mathscr{F}_2)_+\wedge A]^G\cong
\begin{cases}
0,\text{ if }[|k+m|/2]+m\text{ is odd,}\\
[S^{(k+m)+n\lambda},A]^{C_p},\text{ if }[|k+m|/2]+m\text{ is even.}
\end{cases}$$

\textbf{(c)} When $p$ is inverted, we have
$$[S^{k+m\sigma+n\gamma},\widetilde{E\mathscr{F}}_1\wedge A]^G\cong [S^{k+m\sigma},\widetilde{E\mathscr{F}}_1\wedge A]^G\cong [S^{k+m\sigma},\Phi^{C_p}A]^{C_2}.$$

\textbf{(d)} When $2$ is inverted, we have
$$[S^{k+m\sigma+n\gamma},\widetilde{E\mathscr{F}}_2\wedge A]^G\cong [S^{k+n},\widetilde{E\mathscr{F}}_2\wedge A]^G\cong
[S^{k+n},\Phi^{C_2}A]\oplus[S^{k+n},\Phi^GA].$$

Moreover, all isomorphisms above preserve the symmetric monoidal structure.
\end{theorem}

\begin{corollary} \label{hagrog}
In particular, when $A=HA_G$, we have
$$[S^{k+m\sigma+n\gamma},HA_G]^G\cong [S^{(k+n)+(m+n)\sigma},HA_{C_2}^2]^{C_2}\oplus [S^{k+m\sigma},HA_{C_2}]^{C_2}$$
when $p$ is inverted, and
$$[S^{k+m\sigma+n\gamma},HA_G]^G\cong [S^{k+n},H\Z^2]\oplus
\begin{cases}
0,\text{ if }[|k+m|/2]+m\text{ is odd,}\\
[S^{(k+m)+n\lambda},HA_{C_p}]^{C_p},\text{ if }[|k+m|/2]+m\text{ is even.}
\end{cases}$$
when $2$ is inverted.
\end{corollary}

Parts (a) and (b) in \autoref{generald2p} can be proved in the same way as \autoref{z11} and \autoref{z12}.
\medskip

In order to compute the homotopy of $\widetilde{E\mathscr{F}}\wedge A$, for either $\mathscr{F}=\mathscr{F}_1$ or $\mathscr{F}_2$, we need to consider the equivariant homology with coefficients in $N_\mathscr{F}$. We can use the same proof as in \autoref{burnsidehom} to show:

\begin{proposition} \label{nhom}
Let $X$ be a based $G$-CW complex or a $G$-CW spectrum. Then
$$(HN_\mathscr{F})_*X\cong\bigoplus_{H\notin\mathscr{F}}H\Z_*(\Phi^HX/W_GH)$$
where we only choose one $H$ from each conjugacy class of subgroups in $G$.

As a Mackey functor,
$$\underline{HN_\mathscr{F}}_*X(G/L)\cong\bigoplus_{K\notin\mathscr{F}}H\Z_*(\Phi^KX/W_LK)$$
where we only choose one $K$ from each conjugacy class of subgroups in $L$.
\end{proposition}
\bigskip

We will prove:

\begin{proposition} \label{efequiv}
\textbf{(a)} When $\mathscr{F}=\mathscr{F}_1$ and $p$ is inverted, $S^\gamma$ is $\widetilde{E\mathscr{F}}$-equivalent to $S^0$.

\textbf{(b)} When $\mathscr{F}=\mathscr{F}_2$ and $2$ is inverted, $S^\gamma, S^\sigma$ are $\widetilde{E\mathscr{F}}$-equivalent to $S^1,S^0$, respectively.
\end{proposition}

\paragraph{Proof:} Notice that
$$HN_\mathscr{F}=HA_G\wedge\widetilde{E\mathscr{F}}.$$

Since all bounded below $G$-spectra are $HA_G$-local, it suffices to prove \autoref{efequiv} with $\widetilde{E\mathscr{F}}$ replaced by $HN_\mathscr{F}$.
\medskip

\textbf{Part (a):} Assume that $\mathscr{F}=\mathscr{F}_1$ and $p$ is inverted.

\autoref{nhom} shows that:
$$\underline{HN_\mathscr{F}}_*X(G/e)=\underline{HN_\mathscr{F}}_*X(G/C_2)=0,$$
$$\underline{HN_\mathscr{F}}_*X(G/C_p)=H\Z_*(\Phi^{C_p}X),$$
$$\underline{HN_\mathscr{F}}_*X(G/G)=H\Z_*(\Phi^GX)\oplus H\Z_*(\Phi^{C_p}X/G).$$

Since these homology groups only depend on the $C_p$-fixed cells of $X$, $S^0\hookrightarrow S^\gamma$ is an $HN_\mathscr{F}$-equivalence.
\bigskip

\textbf{Part (b):} Assume that $\mathscr{F}=\mathscr{F}_2$ and $2$ is inverted.

\autoref{nhom} shows that:
$$\underline{HN_\mathscr{F}}_*X(G/e)=\underline{HN_\mathscr{F}}_*X(G/C_p)=0,$$
$$\underline{HN_\mathscr{F}}_*X(G/C_2)=H\Z_*(\Phi^{C_2}X),$$
$$\underline{HN_\mathscr{F}}_*X(G/G)=H\Z_*(\Phi^GX)\oplus H\Z_*(\Phi^{C_2}X/G).$$
for any choice of $C_2\subset G$.

These homology groups only depend on the $C_2$-fixed cells of $X$. So the following maps among $D_{2p}$-spaces become $HN_\mathscr{F}$-equivalences:
$$S^\gamma\rightarrow S^\gamma/C_p=S^{1+\sigma}\hookleftarrow S^1,$$
$$S^0\hookrightarrow S^\sigma.$$
$\Box$
\bigskip

\paragraph{Proof of \autoref{generald2p}(c)(d):} Part (c) follows after the fact that $S^{k+m\sigma+n\gamma}$ is $\widetilde{E\mathscr{F}_1}$-equivalent to $S^{k+m\sigma}$ according to \autoref{efequiv}.
\medskip

For part (d), \autoref{efequiv} shows that $S^{k+m\sigma+n\gamma}$ is $\widetilde{E\mathscr{F}_2}$-equivalent to $S^{k+n}$. Thus
$$[S^{k+m\sigma+n\gamma},\widetilde{E\mathscr{F}}_2\wedge A]^G\cong [S^{k+n},\widetilde{E\mathscr{F}}_2\wedge A]^G.$$

Let $\mathscr{F}^\prime$ be the family of all proper subgroups of $G$. Notice that $\widetilde{E\mathscr{F}_2}$ and
$$\widetilde{E\mathscr{F}^\prime}\vee(G/C_2)_+\wedge\widetilde{E(G/C_p)}$$
have the same homotopy types on their fixed point subspaces. \autoref{elmendorf} provides a weak equivalence between these two spaces. Together with \autoref{ringsplit}, we have
$$[S^{k+n},\widetilde{E\mathscr{F}}_2\wedge A]^G\cong[\widetilde{E\mathscr{F}}_2\wedge S^{k+n},A]^G$$
$$\cong[\widetilde{E\mathscr{F}^\prime}\wedge S^{k+n},A]^G\oplus[(G/C_2)_+\wedge\widetilde{E(G/C_p)}\wedge S^{k+n},A]^G$$
$$\cong[\widetilde{E\mathscr{F}^\prime}\wedge S^{k+n},A]^G\oplus[\widetilde{E(G/C_p)}\wedge S^{k+n},A]^{C_2}$$
$$\cong [S^{k+n},\widetilde{E\mathscr{F}^\prime}\wedge A]^G\oplus[S^{k+n},\widetilde{E(G/C_p)}\wedge A]^{C_2}\cong[S^{k+n},\Phi^GA]\oplus[S^{k+n},\Phi^{C_2}A].$$
$\Box$

\section{Splittings for more general finite groups}

In this section, we will generalize \autoref{twosplit} to more general finite groups $G$.
\medskip

Let $G=G_2\ltimes G_1$. In other words, $G_1\lhd G$, $G_2\subset G$, and $G_2\cong G/G_1$.

Let $\mathscr{F}_1$ be the family which contains all subgroups which are sub-conjugate to $G_2$, $\mathscr{F}_2$ be the family which contains all subgroups of $G_1$.

\begin{proposition} \label{gen11}
Assume that $|G_1|,|G_2|$ are relatively prime. 

When $|G_1|$ is inverted, $\underline{H}_*^G(E\mathscr{F}_1;A_G)$ is concentrated in degree $0$.
\end{proposition}

\begin{proposition} \label{gen21}
When $|G_2|$ is inverted, $\underline{H}_*^G(E\mathscr{F}_2;A_G)$ is concentrated in degree $0$.
\end{proposition}

\begin{proposition} \label{gen12}
When $\mathscr{F}=\mathscr{F}_1$, the denominators of all $c_H$ only contain prime factors dividing $|G_1|$.
\end{proposition}

\begin{proposition} \label{gen22}
When $\mathscr{F}=\mathscr{F}_2$, the denominators of all $c_H$ only contain prime factors dividing $|G_2|$.
\end{proposition}

\begin{remark}
We do not require $|G_1|,|G_2|$ to be relatively prime for the last three propositions. We guess this requirement can also be removed for \autoref{gen11}, but do not have a proof at current time.
\end{remark}

According to \autoref{invertprime}, we can generalize \autoref{twosplit}:

\begin{theorem} \label{bigtwosplit}
Assume that $|G_1|,|G_2|$ are relatively prime. The splitting in \autoref{ringsplit} happens when

\textbf{(a)} $\mathscr{F}=\mathscr{F}_1$ and $|G_1|$ is inverted, or

\textbf{(b)} $\mathscr{F}=\mathscr{F}_2$ and $|G_2|$ is inverted.
\end{theorem}

\begin{remark}
When either $|G_1|$ or $|G_2|$ is inverted, it is possible to use the same techniques as in Sections 5 and 6 to decompose computations about $G$-spectra into computations about $G_1$ or $G_2$-spectra. Then we can glue back these two splittings and recover information of $G$-spectra.
\end{remark}
\medskip

We will prove Proposition \ref{gen11} to \ref{gen22} in the rest of the paper from the easiest to the hardest:
$$\ref{gen21}\rightarrow\ref{gen22}\rightarrow\ref{gen11}\rightarrow\ref{gen12}.$$
The proof of \autoref{gen12} will be given in a separate section according to its length. The cases of $\mathscr{F}_1$ and $\mathscr{F}_2$ are quite different since only $G_1$ is a normal subgroup of $G$.
\medskip

\paragraph{Proof of \autoref{gen21}:} Assume that $\mathscr{F}=\mathscr{F}_2$ and $|G_2|$ is inverted.

Recall that $E\mathscr{F}$ is characterized by its fixed point subspaces as $(E\mathscr{F})^H\simeq *$ if $H\subset G_1$, and $(E\mathscr{F})^H=\emptyset$ otherwise. Thus we can construct $E\mathscr{F}$ explicitly as $ E(G/G_1)$.

According to \autoref{burnsidehom}, $\underline{H}_*(E\mathscr{F};A_G)$ is determined by the homology of different orbit spaces of fixed point subspaces. However, the fixed point subspaces of $E(G/G_1)$ are either empty or the whole space. Thus the orbit spaces only have the form $E(G/G_1)/H\simeq BH$ for some $H\subset G/G_1$. According to \autoref{orduni}, we know that the positive-degree homology groups of these orbit spaces vanish when $|G_2|=|G/G_1|$ is inverted. $\Box$
\bigskip

\paragraph{Proof of \autoref{gen22}:} Let $\mathscr{F}=\mathscr{F}_2$. The coefficients $c_H$ are computed by
$$\sum_H c_Hs_{(K,H)}=1, \text{ }\forall K\in\mathscr{F}$$
with one $H$ chosen from each conjugacy class in $\mathscr{F}$. We can check that $c_{G_1}=|G_2|^{-1}$ and $c_H=0$ for all other $H$. Recall that $s_{(K,H)}$ is the product of $|W_GH|$ with the number of subgroups conjugate to $H$ and containing $K$. Thus for any $K\subset G_1$, $s_{(K,G_1)}=|G_2|$ since $G_1$ is normal and $|W_GG_1|=|G/G_1|=|G_2|$. We have
$$\sum_H c_Hs_{(K,H)}=c_{G_1}s_{(K,G_1)}+\sum_{G_1\neq H\in\mathscr{F}}c_Hs_{(K,H)}=|G_2|^{-1}\cdot|G_2|+0=1.$$

\autoref{gen22} is true since the only nontrivial denominator of any $c_H$ is $|G_2|$. $\Box$
\bigskip

For the other two propositions, we will need the following definition:

\begin{definition} \label{commsub}
For any $K\subset G_2$, define $P_K\subset G_1$ to be
$$P_K:=\{g\in G_1:\text{ }gh=hg\textbf{ },\forall h\in K\}.$$
\end{definition}

\begin{lemma} \label{completecomm}
For any $K\subset G_2$ and $g\in G_1$, $gKg^{-1}\subset G_2$ if and only if $g\in P_K$.
\end{lemma}

\paragraph{Proof:} For any $h\in K$, if $ghg^{-1}\in G_2$, then $ghg^{-1}h^{-1}\in G_2h^{-1}=G_2$. However, $g(hg^{-1}h^{-1})\in G_1$ since $G_1\lhd G$. Since $G_1\cap G_2=\{1\}$, we must have $ghg^{-1}h^{-1}=1$ and hence $gh=hg$, $g\in P_K$. $\Box$

\begin{lemma} \label{nnaction}
The subgroup $P_K$ is closed under the action of $N_{G_2}K$ by conjugation.
\end{lemma}

\paragraph{Proof:} For any $g\in P_K$, $h\in N_{G_2}K$, and $k\in K$, we have $h^{-1}kh\in K$. Thus
$$(h^{-1}gh)k=h^{-1}(g(hkh^{-1}))h=h^{-1}((hkh^{-1})g)h=k(h^{-1}gh)$$
and hence $h^{-1}gh\in P_K$. $\Box$
\bigskip

In order to prove \autoref{gen11}, we need a more general version of \autoref{eqorduni}:

\begin{lemma} \label{genuniorb}
Let $X$ be a $G$-CW complex whose cells only have isotropy groups sub-conjugate to $G_2$. Consider the orbit map $X\rightarrow X/G_1$, which is viewed as a map between $G_2$-spaces. Then for any $H\subset G_2$, the induced map
$$H_*(X^H;\Z)\rightarrow H_*((X/G_1)^H;\Z)$$
becomes a surjection when $|G_1|$ is inverted.

Moreover, if all such maps are also injective, then $X\rightarrow(X/G_1)$ is a $G_2$-equivalence.
\end{lemma}

\paragraph{Proof:} First, for any non-equivariant cell $e$ with isotropy group $aKa^{-1}$ for some $K\subset G_2$ and $a\in G$, write $a=gh$ with $g\in G_1$ and $h\in G_2$. Then $aKa^{-1}=g(hKh^{-1})g^{-1}$ and hence $hKh^{-1}$ is the isotropy group of $g^{-1}e$. Therefore, each $G_1$-cell in $X$ contains at least one non-equivariant cell with isotropy group inside $G_2$.
\medskip

Consider any cell $e$ in $X$ with isotropy group $K\subset G_2$. Fix $H\subset G_2$. If the image of $e$ in $X/G_1$ is $H$-fixed, then the action of $H$ is closed inside $\{ge:g\in G_1\}$. 

This means that for any $h\in H$, $g_1\in G_1$, there exists $g_2\in G_1$ such that $hg_1e=g_2e$. So $g_2^{-1}hg_1\in K$, and hence $g_2^{-1}(hg_1h^{-1})\in Kh^{-1}\subset G_2$. Since $hg_1h^{-1}$ is also an element in $G_1$, 
$$g_2^{-1}(hg_1h^{-1})\in G_1\cap G_2=\{1\}.$$

So $g_2=hg_1h^{-1}$. Then $g_2^{-1}hg_1\in K$ becomes $h\in K$. So the image of $e$ in $X/G_1$ is $H$-fixed if and only if $H\subset K$.

Now we assume that $H\subset K$. For any $g\in G_1$, $ge$ is $H$-fixed if and only if $hge=ge$ for all $h\in H$. Thus $g^{-1}hge=e$ and hence $g^{-1}hg\in K$, $g^{-1}Hg\subset K\subset G_2$. \autoref{completecomm} shows that $g\in P_H$. In this case, $g^{-1}Hg=H\subset K$ is always satisfied. So we have
$$\{ge:g\in G_1\}^H=\{ge:g\in P_H\}.$$
\bigskip

Now we can play the same trick as in \autoref{eqorduni}. Define a map between cellular chain complexes
$$C_*((X/G_1)^H;\Z)\rightarrow C_*(X^H;\Z)$$
as follows: For any $G_1$-cell $G/K_+\wedge e$ of $X$ with $K\supset H$, define the map above by $$(G/K_+\wedge e)/G_1\mapsto\sum_{g\in P_H}ge\in C_*(X^H;\Z).$$
Then the composition
$$C_*((X/G_1)^H;\Z)\rightarrow C_*(X^H;\Z)\rightarrow C_*((X/G_1)^H;\Z)$$
becomes the multiplication by $|P_H|$, which can be passed to homology. Thus when $|G_1|$ is inverted, since $P_H\subset G_1$, $|P_H|$ is also inverted. The above composition becomes an isomorphism, and
$$H_*(X^H;\Z)\rightarrow H_*((X/G_1)^H;\Z)$$
becomes a surjection. $\Box$

\begin{remark}
The above argument still works when we replace $G_1, G_2$ by any subgroups $P\subset G_1$, $K\subset G_2$. Thus if $X$ is a $(K\ltimes P)$-CW complex with certain conditions, $X\rightarrow X/P$ will also be a $K$-equivalence when $|P|$ is inverted.
\end{remark}

We also need the powerful Schur–Zassenhaus theorem:

\begin{theorem} \label{schur}
(\textbf{[Iss08]} 3B) Let $K$ be a finite group and $K_1$ be a normal subgroup of $K$ such that $|K_1|$ and $|K/K_1|$ are relative prime. Then there exists $K_2\subset K$, such that $K=K_2\ltimes K_1$. Moreover, all such choices of $K_2$ are conjugate to each other.
\end{theorem}

\paragraph{Proof of \autoref{gen11}:} Assume that $|G_1|$ is inverted and $\mathscr{F}=\mathscr{F}_1$. 

According to \autoref{burnsidehom}, it suffices to show that $H_*(E\mathscr{F}^H/L;\Z)$ is concentrated in degree $0$, for any $H\lhd L\subset G$. Notice that $E\mathscr{F}^H=\emptyset$ if $H\notin\mathscr{F}$. Moreover, for any $g\in G$, $E\mathscr{F}^H/L$ is homotopy equivalent to $E\mathscr{F}^{g^{-1}Hg}/g^{-1}Lg$. Since $\mathscr{F}$ is the family of all subgroups sub-conjugate to $G_2$, we can assume $H\subset G_2$.

List all elements in $L$ as $L=\{a_1b_1,a_2b_2,...,a_nb_n\}$ with each $a_i\in G_1$ and $b_i\in G_2$. Let $\overline{L}$ be the subset of $G_2$ consisting of all $b_i$. Then $\overline{L}$ is closed under multiplication and inverses since
$$(a_ib_i)(a_jb_j)=(a_i(b_ia_jb_i^{-1}))(b_ib_j),$$
$$(a_ib_i)^{-1}=b_i^{-1}a_i^{-1}=(b_i^{-1}a_i^{-1}b_i)b_i^{-1}.$$

Thus $\overline{L}$ is a subgroup of $G_2$. The equations above also show that the map $L\rightarrow \overline{L}$ sending each $a_ib_i$ to $b_i$ is a surjective group homomorphism. Its kernel $L_0$ is the collection of all $a_ib_i$ with $b_i=1$. Thus $L_0$ is a subgroup of $G_1$ and a normal subgroup of $L$. Moreover, for any $i$, since $H\lhd L$, we have $b_i^{-1}a_i^{-1}Ha_ib_i=H$ and hence $a_i^{-1}Ha_i=b_iHb_i^{-1}\subset G_2$. \autoref{completecomm} tells us $a_i\in P_H$ and hence $b_i\in W_{G_2}H$.

Since $W_{G_2}H$ acts on $P_H$ by conjugation according to \autoref{nnaction}, the group generated by $P_H$ and $\overline{L}\subset W_{G_2}H$ is the semi-direct product $\overline{L}\ltimes P_H$.
\smallskip

Notice that $L_0\lhd L$ and $|L_0|$ are relatively prime to $|L/L_0|=|\overline{L}|$ (since $L_0,\overline{L}$ are subgroups of $G_1,G_2$). Apply \autoref{schur} on $L$ and $L_0$. We have $L=L^\prime\ltimes L_0$ for some $L^\prime$ isomorphic to $\overline{L}$. Notice that $L^\prime\subset L\subset\overline{L}\ltimes P_H$. Since $|L^\prime|=|\overline{L}|$ is relatively prime to $P_H$, the elements in $L^\prime$ are in different $P_H$-cosets. Thus $\overline{L}\ltimes P_H=L^\prime\ltimes P_H$. Apply \autoref{schur} again on $\overline{L}\ltimes P_H$, we know that $\overline{L}$ is conjugate to $L^\prime$. So there exists $g\in P_H$ such that $g^{-1}L^\prime g=\overline{L}$.
\smallskip

Now we have $g^{-1}Lg=\overline{L}\ltimes g^{-1}L_0g$ and
$$E\mathscr{F}^H/L\cong E\mathscr{F}^{g^{-1}Hg}/g^{-1}Lg=E\mathscr{F}^H/(\overline{L}\ltimes g^{-1}L_0g).$$
So we can assume $L=\overline{L}\ltimes L_0$ at the beginning. Otherwise we replace it by $\overline{L}\ltimes g^{-1}L_0g$.
\medskip

Recall that all points in $E\mathscr{F}$ have isotropy groups sub-conjugate to $G_2$. Thus when we view $X=E\mathscr{F}^H$ as an $L$-CW complex, all points have isotropy groups sub-conjugate to $\overline{L}$. So all cells are $L_0$-free. Apply \autoref{genuniorb} with $G=G_2\ltimes G_1$ replaced by $L=\overline{L}\ltimes L_0$. The maps between homology groups in \autoref{genuniorb} are always injective since all fixed point subspaces of $E\mathscr{F}^H$ are either contractible or empty. Thus the projection $E\mathscr{F}^H\rightarrow E\mathscr{F}^H/L_0$ is an $\overline{L}$-equivalence. So we have
$$E\mathscr{F}^H/L=(E\mathscr{F}^H/L_0)/\overline{L}\simeq E\mathscr{F}^H/\overline{L}$$

On the other hand, we have a $G_2$-equivalence $E\mathscr{F}\simeq *$ by checking all fixed point subspaces. Since both $H,\overline{L}$ are subgroups of $G_2$, $E\mathscr{F}^H/\overline{L}\simeq *$. Therefore, the homology of $E\mathscr{F}^H/L$ is concentrated in degree 0. $\Box$
\bigskip

\section{Proof of \autoref{gen12}}

Assume that $|G_1|$ is inverted and $\mathscr{F}=\mathscr{F}_1$ in this section.

Recall that $c_H$ is defined for each $H\in\mathscr{F}$ inductively: For any $K\in\mathscr{F}$, we have
$$\sum_H c_Hs_{(K,H)}=1$$
where we choose one $H$ from each conjugacy class inside $\mathscr{F}$, and $s_{(K,H)}=|(G/H)^K|$. We want to show that the denominators of $c_H$ for all $H\subset G$ are invertible.

Since $\mathscr{F}$ is the family of subgroups sub-conjugate to $G_2$, and $c_H=c_{H^\prime}$ when $H$ is conjugate to $H^\prime$, we can consider $K,H$ in the equation above as subgroups of $G_2$.

For any $abH\in G/H$ with $a\in G_1$, $b\in G_2$, if it is $K$-fixed, then $KabH=abH$ and hence
$$(a^{-1}Ka)(bHb^{-1})=(bHb^{-1})\Rightarrow a^{-1}Ka\subset bHb^{-1}\subset G_2.$$
\autoref{completecomm} shows that $a\in P_K$. So $KabH=abH$ becomes $KbH=bH$ and hence $bH\in(G_2/H)^K$. So we have
$$s_{(K,H)}=|P_K|\cdot|(G_2/H)^K|$$
where $|(G_2/H)^K|$ can be expressed as the product between $|W_{G_2}H|$ and the number of subgroups of $G_2$ which are conjugate to $H$ inside $G_2$ and contain $K$.

A direct consequence of \autoref{completecomm} shows that if $H_1,H_2\subset G_2$, then they are conjugate as subgroups of $G$ if and only if they are conjugate as subgroups of $G_2$. So for each $K\subset G_2$, we have
$$1=\sum_{\text{one }H\text{ from each conjugacy class in }G}c_Hs_{(K,H)}=\sum_{\text{one }H\text{ from each conjugacy class in }G_2}c_Hs_{(K,H)}$$
$$=\sum_{\text{one }H\text{ from each conjugacy class in }G_2}c_H\cdot|P_K|\cdot |W_{G_2}H|\cdot\#(\text{subgroups of }G_2\text{ conjugate to }H\text{ containing }K)$$
$$=\sum_{K\subset H\subset G_2}c_H\cdot|P_K|\cdot|W_{G_2}H|.$$
In the last sum, we choose all $H$ with $K\subset H\subset G_2$, not up to conjugacy. We highlight this new relation:

\begin{proposition} \label{chrelation}
For any $K\subset G_2$, we have
$$\sum_{K\subset H\subset G_2}c_H\cdot|P_K|\cdot|W_{G_2}H|=1.$$
Moreover, when $H,H^\prime$ are conjugate subgroups of $G_2$, we have $c_H=c_{H^\prime}$.
\end{proposition}

Before we start computing the value of each $c_H$, we need to study relations between subgroups of $G_1$ and subgroups of $G_2$.

\paragraph{Notation:} We will use letters $K,H,L$ to denote subgroups of $G_2$, and $P,Q$ to denote subgroups of $G_1$ in the rest of this section.
\medskip

First we extend \autoref{commsub}:

\begin{definition} \label{decent}
For any $K\subset G_2$, let
$$P_K:=\{g\in G_1:\text{ }gh=hg,\forall h\in K\}.$$

For any $P\subset G_1$, let
$$K_P:=\{h\in G_2:\text{ }gh=hg,\forall g\in P\}.$$

For any subgroup $K$ of $G_2$, it is called \textbf{G$_1$-good} if $K=K_P$ for some $P\subset G_1$.

For any subgroup $P$ of $G_1$, it is called \textbf{G$_2$-good} if $P=P_K$ for some $K\subset G_2$.
\end{definition}

Some of their properties are listed below:

\begin{proposition} \label{decentprop}
We have:

\textbf{(a)} $P\subset P_{K_P}$, with equality when $P$ is $G_2$-good.

\textbf{(b)} $K\subset K_{P_K}$, with equality when $K$ is $G_1$-good.

\textbf{(c)} For any $K_1,K_2\subset G_2$, $K_1\subset K_2$ implies $P_{K_2}\subset P_{K_1}$. The converse is true if $K_1,K_2$ are $G_1$-good.

\textbf{(d)} For any $P_1,P_2\subset G_1$, $P_1\subset P_2$ implies $K_{P_2}\subset K_{P_1}$. The converse is true if $P_1,P_2$ are $G_2$-good.

\textbf{(e)} For any $b\in G_2$, $K\subset G_2$, $P_{b^{-1}Kb}=b^{-1}P_Kb$. Thus if $K$ is $G_1$-good, so is $b^{-1}Kb$.

\textbf{(f)} If both $P_1,P_2\subset G_1$ are $G_2$-good, so is $P_1\cap P_2$.
\end{proposition}

\paragraph{Proof:} Parts (a) to (e) follow immediately after \autoref{decent}.

For part (f), let $K_1=K_{P_1}$ and $K_2=K_{P_2}$. Parts (a) and (b) show that $P_1=P_{K_1}$ and $P_2=P_{K_2}$.

Let $K$ be the smallest subgroup of $G_2$ which contains both $K_1$ and $K_2$. In other words, it is generated by elements in $K_1$ and $K_2$. We compare $P_K$ and $P_1\cap P_2$:

For any $g\in P_1\cap P_2$, we have $gh=hg$ for all $h\in K_1$ and all $h\in K_2$. Since $K$ is generated by $K_1,K_2$, $gh=hg$ for all $h\in K$. Thus $P_1\cap P_2\subset P_K$.

On the other hand, for any $g\in P_K$, since $K_1,K_2\subset K$, we have $gh=hg$ for all $h\in K_1$ and all $h\in K_2$. Thus $g\in P_1$ and $g\in P_2$. So we have $P_K\subset P_1\cap P_2$. Therefore, $P_1\cap P_2=P_K$ is also $G_2$-good. $\Box$
\bigskip

It suffices to consider $G_1$ and $G_2$-good subgroups since:

\begin{proposition} \label{onlydecent}
Let $\mathbb{D}$ be the collection of all $G_1$-good subgroups of $G_2$. Then $c_H=0$ for any $H\notin\mathbb{D}$.
\end{proposition}

\paragraph{Proof:} We use induction. Clearly $G_2$ itself is $G_1$-good since $G_2=K_{\{1\}}$.

Fix $K\subset G_2$ which is not $G_1$-good. Assume that the proposition is true for all $H\subset G_2$ and $H\supsetneqq K$. Let $K_0=K_{P_K}$. Then $K\subsetneqq K_0$ and $|P_K|=|P_{K_0}|$. \autoref{chrelation} gives us:
$$\sum_{K\subset H\subset G_2}c_H\cdot|W_{G_2}H|=|P_K|^{-1},$$
$$\sum_{K_0\subset H\subset G_2}c_H\cdot|W_{G_2}H|=|P_{K_0}|^{-1}=|P_K|^{-1}.$$
Taking the difference gives
$$\sum_{K\subset H,K_0\nsubseteq H, H\subset G_2}c_H\cdot|W_{G_2}H|=0.$$

For any $H$ which contains $K$ but not $K_0$, \autoref{decentprop}(c) shows that $P_H\subset P_K=P_{K_0}$. \autoref{decentprop}(d) then shows that $H$ cannot be $G_1$-good, since otherwise we must have $K_0\subset H$. Thus all summands in
$$\sum_{K\subset H,K_0\nsubseteq H, H\subset G_2}c_H\cdot|W_{G_2}H|=0$$
have $H\notin\mathbb{D}$. We already know that $c_H=0$ for all $K\subsetneqq H\notin\mathbb{D}$. Thus $c_K$ must also be zero. $\Box$
\bigskip

Now we can refine \autoref{chrelation}:

\begin{proposition} \label{finalrelation}
For any $K\in\mathbb{D}$, we have
$$\sum_{K\subset H\in\mathbb{D}}c_H\cdot|P_K|\cdot|W_{G_2}H|=1.$$
Moreover, $c_H=c_{H^\prime}$ if $H$ is conjugate to $H^\prime$. If $H\notin\mathbb{D}$, $c_H=0$.
\end{proposition}
\bigskip

We have finished all preparations and can start proving \autoref{gen12}.
\bigskip

We can still use induction. For $K=G_2$, \autoref{finalrelation} shows that
$$c_{G_2}\cdot|P_{G_2}|=1.$$
Since $P_{G_2}\subset G_1$, $|P_{G_2}|$ becomes invertible after we invert $|G_1|$. Thus the denominator of $c_{G_2}$ is invertible.
\medskip

Now fix $K\subset G_2$. Assume that the denominators of all $c_H$ with $K\subsetneqq H\subset G_2$ are invertible. We will prove the same thing for $c_K$. We can assume $K\in\mathbb{D}$ (otherwise $c_K=0$ by \autoref{finalrelation}).
\medskip

Let $P=P_K$. List all maximal proper $G_2$-good subgroups of $P$ as $P_1,P_2,...,P_n$. For any $S\subset\{1,2,...,n\}$. Define
$$P_S:=P\cap\left(\bigcap_{i\in S}P_i\right).$$
\autoref{decentprop}(f) shows that $P_S$ is $G_2$-good. Define $K_S:=K_{P_S}$. Then we have
$$P_{\{i\}}=P_i,\text{ }P_\emptyset=P,\text{ }K_\emptyset=K.$$
\medskip

We can also extend \autoref{nnaction}:

\begin{lemma} \label{naction}
Consider the action of $N_{G_2}K$ on $P_K$ by conjugation. The action of any $h\in N_{G_2}K$ permutes different $P_S$.
\end{lemma}

\paragraph{Proof:} For any $P_i$ and $h\in N_{G_2}K$, notice that $h^{-1}P_ih$ is also a maximal proper $G_2$-good subgroup of $P=P_K$. Thus $h^{-1}P_ih=P_j$ for some $j$. So the action of $h$ becomes a permutation among $P_1,P_2,...,P_n$, and also among general $P_S$, since each $P_S$ is defined as the intersection of some $P_i$. $\Box$
\bigskip

Write $K_i:=K_{\{i\}}=K_{P_i}$, $i=1,2,...,n$. A quick application of \autoref{decentprop} shows:

\begin{lemma} \label{maximalinclude}
For any $K\subsetneqq H\in\mathbb{D}$, $K_i\subset H$ for at least one $i$.
\end{lemma}

Apply \autoref{finalrelation} for each $S\subset\{1,2,...,n\}$:
$$\sum_{K_S\subset H\in\mathbb{D}}c_H\cdot|P_S|\cdot|W_{G_2}H|=1.$$
Add these equations together with additional signs:
$$\sum_{S}(-1)^{|S|}\left(\sum_{K_S\subset H\in\mathbb{D}}c_H\cdot|P_S|\cdot|W_{G_2}H|\right)=\sum_S(-1)^{|S|}=\binom{n}{0}-\binom{n}{1}+...+(-1)^n\binom{n}{n}=0.$$

Rewriting the sum in the left hand side:
$$\sum_{K\subset H\in\mathbb{D}}c_H\cdot|W_{G_2}H|\cdot\left(\sum_{K_S\subset H}(-1)^{|S|}|P_S|\right)=0.$$

Separate the summand of $K$:
$$-c_K\cdot|W_{G_2}K|\cdot|P_K|=\sum_{K\subsetneqq H\in\mathbb{D}}c_H\cdot|W_{G_2}H|\cdot\left(\sum_{K_S\subset H}(-1)^{|S|}|P_S|\right).$$

We divide all $H$ with $K\subsetneqq H\in\mathbb{D}$ into different classes: $H_1,H_2$ are in the same class if there exists $h\in N_{G_2}K$ such that $H_2=h^{-1}H_1h$.

For $K\subsetneqq H\in\mathbb{D}$ and $h\in N_{G_2}K$, $h^{-1}Hh$ is also a $G_1$-good subgroup which has $K$ as a proper subgroup. Moreover, for any $h_1,h_2\in N_{G_2}K$, we have $h_1^{-1}Hh_1=h_2^{-1}Hh_2$ if and only if $h_1,h_2$ are in the same $N_{G_2}H$-coset. Thus the class of $H$ contains $|N_{G_2}K|/|N_{G_2}K\cap N_{G_2}H|$ elements.

For any $H_1,H_2$ in the same class, since $H_1$ is conjugate to $H_2$, we have
$$c_{H_1}=c_{H_2},\text{ }|W_{G_2}H_1|=|W_{G_2}H_2|.$$
\autoref{naction} shows that
$$\sum_{K_S\subset H_1}(-1)^S|P_S|=\sum_{K_S\subset H_2}(-1)^S|P_S|.$$ 
Thus we have
$$-c_K\cdot|W_{G_2}K|\cdot|P_K|=\sum_{K\subsetneqq H\in\mathbb{D}}c_H\cdot|W_{G_2}H|\cdot\left(\sum_{K_S\subset H}(-1)^{|S|}|P_S|\right)$$
$$=\sum_{\text{one }H\text{ from each class}}c_H\cdot|W_{G_2}H|\cdot(\text{size of the class})\cdot\left(\sum_{K_S\subset H}(-1)^{|S|}|P_S|\right)$$
$$=\sum_{\text{one }H\text{ from each class}}c_H\cdot\frac{|N_{G_2}H|}{|H|}\cdot\frac{|N_{G_2}K|}{|N_{G_2}K\cap N_{G_2}H|}\cdot\left(\sum_{K_S\subset H}(-1)^{|S|}|P_S|\right).$$

When $|G_1|$ is inverted, $|P_K|$ becomes invertible. We also know that $c_H$ has invertible denominators for all $H\supsetneqq K$ by assumption. Now it suffices to prove:

\begin{proposition} \label{divide}
For any $K\subsetneqq H\in\mathbb{D}$, $|W_{G_2}K|$ divides 
$$\frac{|N_{G_2}H|}{|H|}\cdot\frac{|N_{G_2}K|}{|N_{G_2}K\cap N_{G_2}H|}\cdot\left(\sum_{K_S\subset H}(-1)^{|S|}|P_S|\right).$$
\end{proposition}

\paragraph{Proof:} Notice that $H$ must contains at least one $K_i$ according to \autoref{maximalinclude}. Without loss of generality, assume that $H$ contains $K_1,K_2,...,K_m$, but does not contain $K_{m+1},K_{m+2},...,K_n$. Then \autoref{decentprop} shows that $P_H$ is contained in $P_1,P_2,...,P_m$, but not in $P_{m+1},...,P_n$.

We have
$$\sum_{K_S\subset H}(-1)^{|S|}|P_S|=\sum_{S\subset\{1,2,...,m\}}(-1)^{|S|}|P_S|=|P|+\sum_{\emptyset\neq S\subset\{1,2,...,m\}}(-1)^{|S|}\left|\bigcap_{i\in S}P_i\right|.$$
This is the number of elements in $P$ which are not contained in any of $P_1,P_2,...,P_m$.
\medskip

Now we consider the group action on $P_K$. But instead of the whole $N_{G_2}K$, we only use the action of $L:=N_{G_2}K\cap H$.

For any $g\in P-\bigcup_{i=1}^mP_i$ and $h\in L$, assume that $h^{-1}gh\in P_i$ for some $i\in\{1,2,...,m\}$. Then $g\in hP_ih^{-1}$. However, since $h\in N_{G_2}K$, $hP_ih^{-1}$ is another maximal proper $G_2$-good subgroup of $P$, denoted by $P_j$. Since $h\in H$ and $K_i\subset H$, $hK_ih^{-1}\subset H$ and hence $K_j\subset H$. This means that $j\in\{1,2,...,m\}$ and $g\in P_j$, which is a contradiction.

In other words, the subset $P-\bigcup_{i=1}^mP_i$ is closed under the action of $L$.
\medskip

We can also describe the orbit types of this $L$-action.
\smallskip

For any $g\in P-\bigcup_{i=1}^mP_i$ and $h\in L$, assume that $h^{-1}gh=g$ and $h\notin K$. Then $g\in P_{K^\prime}\subsetneqq P$, where $K^\prime$ is generated by $K$ and $h$. Since $h\in H$ and $K\subset H$, we have $K^\prime\subset H$ and hence $P_H\subset P_{K^\prime}\subsetneqq P_K$. Since $P_1,P_2,...,P_m$ are all the maximal $G_2$-good strict subgroups of $P=P_K$ which contain $P_H$, $P_{K^\prime}$ must be contained in one of them. Thus $g$ is contained in one of $P_1,...,P_m$, which is a contradiction.

Therefore, $h^{-1}gh=g$ if and only if $h\in K$. Thus $P-\bigcup_{i=1}^mP_i$ only have $L$-orbits as $L/K$. So we have
$$|L/K|=\frac{|N_{G_2}K\cap H|}{|K|}\text{ divides }\left|P-\bigcup_{i=1}^mP_i\right|=\sum_{K_S\subset H}(-1)^{|S|}|P_S|.$$
\smallskip

Thus it suffices to show that
$$|W_{G_2}K|=\frac{|N_{G_2}K|}{|K|}\text{ divides }\frac{|N_{G_2}H|}{|H|}\cdot\frac{|N_{G_2}K|}{|N_{G_2}K\cap N_{G_2}H|}\cdot\frac{|N_{G_2}K\cap H|}{|K|},$$
or equivalently,
$$\frac{|N_{G_2}K\cap N_{G_2}H|}{|N_{G_2}K\cap H|}\text{ divides }\frac{|N_{G_2}H|}{|H|}=|(N_{G_2}H)/H|.$$

This relation follows from the fact that $N_{G_2}K\cap H$ is the kernel of the group homomorphism
$$N_{G_2}K\cap N_{G_2}H\hookrightarrow N_{G_2}H\rightarrow (N_{G_2}H)/H.$$
$\Box$

\bibliographystyle{abbrv}

\begin{thebibliography}{13}

\bibitem[Bar08]{B}{D. J. Barnes, \textit{Rational Equivariant Spectra}, PhD Thesis.}

\bibitem[Dieck72]{D}{T. tom Dieck, \textit{Orbittypen und äquivariante Homologie}, I. Arch. Math. 23, pp 307-317, 1972.}

\bibitem[Elm83]{E}{A. D. Elmendorf, \textit{Systems of fixed point sets}, Transactions of the AMS, Volume 277, No. 1, 1983.}

\bibitem[GM95]{G}{J. P. C. Greenlees, J. P. May, \textit{Generalized Tate cohomology}, Memoirs of the AMS, 1995.}

\bibitem[HHR17]{H}{M. A. Hill, M. J. Hopkins, D. C. Ravenel, \textit{The slice spectral sequence for the $C_4$ analog of real $K$-theory}, Forum Math., 29(2), pp 383-447, 2017.}

\bibitem[Iss08]{I}{I. M. Issac, \textit{Finite group theory}, Graduate Studies in Mathmetics, Volume 92, 2008.}

\bibitem[KL20]{K}{I. Kriz, Y. Lu, \textit{On the $RO(G)$-graded coefficients of dihedral equivariant cohomology}, arXiv:2005.01225v1.}

\bibitem[LMS86]{L}{L. G. Lewis, J. P. May, M. Steinberger, \textit{Equivariant stable homotopy theory}, Springer-Verlag, 1986.}

\bibitem[LM04]{L}{L. G. Lewis, M. A. Mandell, \textit{Equivariant universal coefficient and Künneth spectral sequences}, arXiv:math/0410162v1.}

\bibitem[May96]{M}{J. P. May, \textit{Equivariant homotopy and cohomology theory}, Conference Board of the Mathematical Sciences, No. 91, 1996.}

\bibitem[Wil75]{W}{S. J. Willson, \textit{Equivariant homology theories on $G$-complexes}, Trans. Amer. Math. Soc. 212, 1975, pp 155-271.}

\bibitem[Zen18]{Z}{M. Zeng, \textit{Mackey functors, equivariant Eilenberg-Maclane spectra and their slices}, PhD thesis.}

\bibitem[Zou18]{Z}{Y. Zou, \textit{$RO(D_{2p})$-graded Slice Spectral Sequence of $H\underline{\Z}$}, PhD Thesis.}


\end{thebibliography}

\end{document}